\documentclass{amsart}
\usepackage{amsmath}
\usepackage{amsthm}
\usepackage{amssymb}

 \newtheorem{thm}{Theorem}[subsection]
 \newtheorem{cor}[thm]{Corollary}
 \newtheorem{lem}[thm]{Lemma}
 \newtheorem{prop}[thm]{Proposition}
 \theoremstyle{definition}
 \newtheorem{defn}[thm]{Definition}
 \theoremstyle{remark}
 \newtheorem{rem}[thm]{Remark}
 \theoremstyle{example}
 
 \numberwithin{equation}{subsection}
\newcommand{\pn}{\noindent}

\newcommand{\ZZ}{\mathbb{Z}}
\newcommand{\QQ}{\mathbb{Q}}
\newcommand{\RR}{\mathbb{R}}
\newcommand{\LL}{\mathbb{L}}
\newcommand{\CC}{\mathbb{C}}
\newcommand{\GG}{\mathbb{G}}

\newcommand{\li}{\mathrm{Lie}\,}

\newcommand{\Hom}{\mathrm{Hom}}
\newcommand{\Ext}{\mathrm{Ext}}
\newcommand{\Biext}{\mathrm{Biext}}
\newcommand{\uHom}{\underline{\mathrm{Hom}}}
\newcommand{\uExt}{\underline{\mathrm{Ext}}}

\newcommand{\bBiext}{\mathrm{\mathbf{Biext}}}

\newcommand{\W}{\mathrm{W}}
\newcommand{\F}{\mathrm{F}}
\newcommand{\T}{\mathrm{T}}
\newcommand{\h}{\mathrm{H}}

\newcommand{\Gr}{\mathrm{Gr}}
\newcommand{\rk}{\mathrm{rk}}
\newcommand{\an}{\mathrm{an}}
\newcommand{\spec}{{\mathrm{Spec}}\,}
\newcommand{\ok}{\overline k}

\newcommand{\gal}{{\mathrm{Gal}}(\ok/k)}

\newcommand{\Sp}{{\mathrm{Sp}}\,}
\newcommand{\dR}{\mathrm{dR}}

\begin{document}

\title[Multilinear morphisms]
{Multilinear morphisms between 1-motives}

\author{Cristiana Bertolin}

\address{NWF-I Mathematik, Universit\"at Regensburg, D-93040 Regensburg}

\email{cristiana.bertolin@mathematik.uni-regensburg.de}



\date{}
\dedicatory{}

\commby{Cristiana Bertolin}


\begin{abstract}

We introduce the notion of biextensions of 1-motives over an arbitrary scheme $S$ and we define bilinear morphisms between 1-motives as isomorphism classes of such biextensions.
If $S$ is the spectrum of a field
of characteristic 0, we check that these biextensions define bilinear morphisms between the realizations of 1-motives. Generalizing we obtain the notion of multilinear morphisms between 1-motives.

\end{abstract}

\maketitle


\section*{Introduction}

Let $S$ be a scheme.
A \textbf{1-motive} $M=(X,A,Y(1),G,u)$ over $S$ consists of
\begin{itemize}
    \item  an $S$-group scheme $X$ which is locally for the \'etale
topology a constant group scheme defined by a finitely generated free
$\ZZ \,$-module,
    \item an extension $G$ of an abelian $S$-scheme $A$ by an $S$-torus $Y(1),$
with co\-cha\-ra\-cter group $Y$,
    \item a morphism $u:X \rightarrow G$ of $S$-group schemes.
\end{itemize}
If $S$ is the spectrum of the field ${\CC}$ of complex numbers, in \cite{D1} (10.1.3) Deligne proves that the category of 1-motives over $S$ is equivalent through the functor ``Hodge realization'' $ M \mapsto {\T}_{\h}(M)$  to the category of $\QQ\,$-mixed Hodge structures $\h$, endowed with a torsion-free $\ZZ \,$-lattice, of type $\{ (0,0),(-1,0),(0,-1), (-1,-1)\}$, and with the quotient ${\rm Gr}_{-1}^{{\W}}(\h)$ polarizable. In the category $\mathcal {MHS}$ of mixed Hodge structures there is an obvious notion of tensor product. If $M_i$ (for $i=1, \ldots, n$) and $M$ are 1-motives defined over $\CC$, the group ${\Hom}_{\mathcal{MHS}}(\otimes_{i=1}^n{\T}_{\h}(M_i) , {\T}_{\h}(M) )$
is hence defined. Our aim in this paper is to show that this group admits a purely algebraic interpretation. More precisely,
if $M_i$ (for $i=1, \ldots, n$) and $M$ are 1-motives defined over an arbitrary scheme $S$, using biextensions we define a group
 $${\Hom}(M_1, \ldots, M_n;M)$$
  of \emph{multilinear morphisms from
$M_1 \times \ldots \times M_n$ to $M$}, which for $S={\spec}(\CC)$ can be identified with the group ${\Hom}_{\mathcal{MHS}}(\otimes_{i=1}^n{\T}_{\h}(M_i) , {\T}_{\h}(M) )$.

One hopes that for any field $k$, there is a $\QQ \,$-linear Tannakian category of mixed motives
over $k$. The category of 1-motives over $k$, taken up to isogeny (i.e. tensorizing the Hom-groups by ${\QQ}$), should be a subcategory, and our notion of multilinear morphisms between 1-motives should agree
with the notion of multilinear morphisms in this Tannakian category.
For 1-motives, we are able to define multilinear morphisms between 1-motives
working integrally and over an arbitrary base scheme $S$ and we check that if $S$ is the spectrum of the field $k$ of characteristic 0 embeddable in $\CC$, our definition agrees with the notion of multilinear morphisms in the integral version of the Tannakian category
${\mathcal{MR}}_{\ZZ}(k)$ of mixed realizations over $k$ introduced by Jannsen in~\cite{J} I 2.1. Our results might give some guidance as to what
to hope for more general mixed motives.

The idea of defining morphisms through biextensions goes back to Gro\-then\-dieck, who
 defines $\ell$-adic pairings from biextensions (cf.~\cite{SGA7} Expos\'e VIII):
if $P,Q,G$ are three abelian groups of a topos \textbf{T}, to each isomorphism class of biextensions of $(P,Q)$ by $G,$ he associates a pairing $({}_{l^n}P)_{n \geq 0} \otimes ({}_{l^n}Q)_{n \geq 0}  \rightarrow ({}_{l^n}G)_{n \geq 0} $ where $({}_{l^n}P)_{n \geq 0} $ (resp. $({}_{l^n}Q)_{n \geq 0} $, $({}_{l^n}G)_{n \geq 0} $) is the projective system constructed from the kernels ${}_{l^n}P$ ( resp.
${}_{l^n}Q$ , ${}_{l^n}G$) of the multiplication by $l^n$ for each $n \geq 0$. Let $K_i=[A_i \stackrel{u_i}{\longrightarrow} B_i]$ (for $i=1,2$) be two complexes of abelian sheaves (over a topos \textbf{T}) concentrated in degree 0 and -1.
Generalizing Grothendieck's work, in~\cite{D1} (10.2.1) Deligne defines the notion of biextension of $(K_1,K_2)$ by an abelian sheaf. Applying this definition to two 1-motives $M_1,M_2$ defined over $\CC$ and to ${\GG}_m,$
he associates, to each isomorphism class of such biextensions, a morphism from the tensor product of the Hodge realizations (resp. the De Rham realizations, resp. the $\ell$-adic realizations) of $M_1$ and $M_2$ to the Hodge realization (resp. the De Rham realization, resp. $\ell$-adic realization) of ${\GG}_m.$

Let $K_i=[A_i \stackrel{u_i}{\longrightarrow} B_i]$ (for $i=1,2,3$) be three complexes of abelian sheaves (over a topos \textbf{T}) concentrated in degree 0 and -1. In this paper we define the notion of biextension of $(K_1,K_2)$ by $K_3$ (see Definition~\ref{defbiext}). In the special case where $A_3=0$, i.e. $K_3=[0 \longrightarrow B_3]$, our definition coincides with Deligne'one ~\cite{D1} (10.2.1). Since we can view 1-motives as
complexes of commutative $S$-group schemes concentrated in degree 0 and -1, applying our definition of biextension of complexes of abelian sheaves concentrated in degree 0 and -1 to 1-motives, we get the following notion of biextension of 1-motives by 1-motives:

\begin{defn}\label{biext1-mot}
 Let $M_i=[X_i \stackrel{u_i}{\longrightarrow}G_i]$ (for $i=1,2,3$) be a 1-motive over a scheme $S$. A biextension $({\mathcal{B}}, \Psi_1, \Psi_2,\lambda)$ of $(M_1,M_2)$ by $M_3$ consists of
 \begin{enumerate}
    \item a biextension of $\mathcal{B}$ of $(G_1, G_2)$ by $G_3$;
    \item a trivialization $\Psi_1 $ (resp. $\Psi_2$) of the biextension $(u_1, id_{G_2})^* {\mathcal{B}}$ of $(X_1, G_2)$ by $G_3$ (resp. of the biextension $(id_{G_1},u_2)^*{\mathcal{B}}$ of $(G_1, X_2)$
by $G_3$) obtained as pull-back of ${\mathcal{B}}$ via $(u_1,id_{G_2})$ (resp. via $(id_{G_1},u_2)$). These two trivializations $\Psi_1 $ and $\Psi_2 $ have to coincide over $(X_1, X_2)$, i.e.
\[(u_1,id_{X_2})^*\Psi_2=\Psi=(id_{X_1},u_2)^*\Psi_1 \]
\pn with $\Psi$ a trivialization of the biextension $(u_1, u_2)^*{\mathcal{B}}$
 of $(X_1, X_2)$ by $G_3$ obtained as pull-back via $(u_1,u_2)$ of the biextension ${\mathcal{B}}$;
    \item a morphism $\lambda:X_1 \otimes X_2 \rightarrow X_3$ of $S$-group schemes such that $u_3 \circ \lambda: X_1 \otimes X_2 \rightarrow G_3$ is compatible with the trivialization  $\Psi$ of the biextension  $(u_1,u_2)^* {\mathcal{B}}$ of $(X_1, X_2)$ by $G_3$.
\end{enumerate}
\end{defn}

We denote by ${\Biext}^1(M_1,M_2;M_3)$ the group of isomorphism classes of biextensions of $(M_1,M_2)$ by $M_3$.

\begin{defn}
Let $M_i$ (for $i=1,2,3$) be a 1-motive over $S$. A morphism
\[ M_1 \otimes M_2 \longrightarrow M_3\]
 from the tensor product of $M_1$ and $M_2$ to the 1-motive $M_3$ is an isomorphism class of biextensions of $(M_1,M_2)$ by $M_3$. Moreover, to $M_1,M_2$ and $M_3$ we associate a group ${\Hom}(M_1, M_2;M_3)$
 defined in the following way:
\[ {\Hom}(M_1, M_2;M_3) := {\Biext}^1(M_1,M_2;M_3), \]
i.e. ${\Hom}(M_1, M_2;M_3)$ is the group of bilinear morphisms from
$ M_1 \times M_2$ to $M_3$.
 \end{defn}

 Observe that the tensor product $M_1 \otimes M_2$ of two 1-motives is not defined yet, and that according to the compatibility between the tensor product of motives and the weight filtration of motives, such a tensor product $M_1 \otimes M_2$ is no longer a 1-motive. But since morphisms of motives have to respect the weight filtration ${\W}_*$, the only non trivial components of the morphism $M_1 \otimes M_2 \rightarrow M_3$ are the components of the morphism
from the 1-motive underlying the quotient $ M_1 \otimes M_2 / {\W}_{-3}(M_1 \otimes M_2)$ to the 1-motive $M_3$. Therefore for our goal only the 1-motive underlying $ M_1 \otimes M_2 / {\W}_{-3}(M_1 \otimes M_2)$ is involved. We construct explicitly this 1-motive in section 2.

Imposing the fact that morphisms of motives have to respect the weight filtration ${\W}_*$,
if $M_i$ (for $i=1,\ldots,n$) and $M$ are 1-motives over $S$, we describe \emph{the group
\[ {\Hom}(M_1, \ldots, M_n;M)  \]
of multilinear morphisms from $M_1 \times \ldots \times M_n$ to $M$} always in terms of biextensions of 1-motives by 1-motives (Theorem~\ref{thmotimes}).

We finish studying the cases in which we can describe the group of isomorphism classes of biextensions of 1-motives as a group of bilinear morphisms in an appropriate category:
\begin{itemize}
  \item If $S={\spec}(\CC)$, the group of isomorphism classes of biextensions of $(M_1,M_2)$ by $M_3$ is isomorphic to
 the group of morphisms (of the category $\mathcal {MHS}$ of mixed Hodge structures) from the tensor product ${\T}_{\h}(M_1) \otimes {\T}_{\h}(M_2)$ of the Hodge realizations of $M_1$ and $M_2$ to the Hodge realization ${\T}_{\h}(M_3)$ of $M_3$:
\begin{equation}\label{intro1}
    {\Hom}(M_1,M_2;M_3) \cong {\Hom}_{\mathcal{MHS}}({\T}_{\h}(M_1) \otimes {\T}_{\h}(M_2) , {\T}_{\h}(M_3) ).
\end{equation}
  \item If $S$ is the spectrum of a field $k$ of characteristic 0 embeddable in $\CC$, modulo isogenies the group of isomorphism classes of biextensions of $(M_1,M_2)$ by $M_3$ is isomorphic to
 the group of morphisms of the category ${\mathcal{MR}}_{\ZZ}(k)$ of mixed realizations over $k$ with integral structure (integral version of the Tannakian category of mixed realizations over $k$ introduced by Jannsen in~\cite{J} I 2.1) from the tensor product ${\T}(M_1)\otimes {\T}(M_2)$ of the realizations of $M_1$ and $M_2$ to the realization ${\T}(M_3)$ of $M_3$:
\begin{equation}\label{intro2}
    {\Hom}(M_1, M_2;M_3) \otimes {\QQ} \cong {\Hom}_{{\mathcal{MR}}_{\ZZ}(k)}\big(
{\T}(M_1)\otimes {\T}(M_2), {\T}(M_3)\big).
\end{equation}
In other words, following Deligne's philosophy of motives described in ~\cite{D3} 1.11,
the notion of biextensions of 1-motives by 1-motives that we have introduced furnishes \emph{the geometrical origin} of the morphisms of ${\mathcal{MR}}_{\ZZ}(k)$
from the tensor product of the realizations of two 1-motives to the realization of another 1-motive, which are therefore \emph{motivic morphisms}.
\end{itemize}

We expect to have a description of biextensions of 1-motives by 1-motives as bilinear morphisms also in the following categories:
\begin{itemize}
  \item If $S$ is a scheme of finite type over $\CC$, we expect to generalize~(\ref{intro1}) finding a description of the group ${\Hom}(M_1, M_2;M_3)$ in terms of bilinear morphisms of an appropriate subcategory of the category of variations of mixed Hodge structures.
  \item If $S$ is a scheme of finite type over $\QQ$,
we expect to generalize~(\ref{intro2}) getting a description of
 the group ${\Hom}(M_1, M_2;M_3) \otimes {\QQ}$ as a group of bilinear morphisms in the Tannakian category $\mathcal{M}(S)$ of mixed realizations over $S$ with integral structure introduced by Deligne in~\cite{D3} 1.21, 1.23 and 1.24. Taking the inductive limit, it should be possible to generalize this last case to any scheme $S$ of characteristic 0.
 \item If $S$ is the spectrum of a perfect field $k$, we expect to get a description of the group ${\Hom}(M_1, M_2;M_3) \otimes {\QQ}$ in terms of bilinear morphisms of the Voevodsky triangulated category $\mathrm{DM}^{\mathrm{eff}}_{\mathrm{gm}}(k)$ of effective geometrical motives, using the Orgogozo-Voevodsky functor from
     the derived category of the category of $1$-motives up to isogeny to the category $\mathrm{DM}^{\mathrm{eff}}_{\mathrm{gm}}(k) \otimes {\QQ}$ (see~\cite{O}).
\end{itemize}

Seeing biextensions as multilinear morphisms was already used in the computation of the unipotent radical of the Lie algebra of the motivic Galois group of a 1-motive defined over a field $k$ of characteristic 0.
In fact in~\cite{B} (1.3.1), using Deligne's definition of biextension
of 1-motives by ${\GG}_m$, we defined a morphism from the tensor product $M_1 \otimes M_2$ of two 1-motives to a torus as an isomorphism class of biextensions of $(M_1,M_2)$ by this torus.

Remark that the results obtained in~\cite{B3}, in particular \textbf{Theorem A}, \textbf{Theorem B} and \textbf{Theorem C}, mean that \emph{biextensions respect the weight filtration ${\W}_*$ of motives}, i.e. they satisfy the main property of morphisms of motives.


\section*{Acknowledgment}

Je tiens \`a remercier Pierre Deligne pour ses commentaires \`a propos de ce travail. Ich bedanke mich bei Uwe Jannsen f\"ur
die sehr interessante und motivierende Diskussionen, die wir w\"ahrend meines Besuches in Regensburg gef\"uhrt haben.

\section*{Notation}

In this paper $S$ is a scheme.

Because of the similar behavior of the different cohomology theories,
it is expected that motives satisfy the following basic properties:
\begin{enumerate}
    \item \emph{The weight filtration ${\W}_*$}: each motive $M$ is endowed with an increasing filtration ${\W}_*$, called the weight filtration. A motive $M$ is said to be \emph{pure of weight $i$} if ${\W}_i(M)=M$ and $W_{i-1}(M)=0.$ Motives which are not pure are called \emph{mixed motives}. This weight filtration ${\W}_*$ is strictly compatible with any morphism $f:M \rightarrow N$ between motives, i.e.
\[ f(M) \cap {\W}_i(N) =  f({\W}_i(M)). \]
In terms of pure motives, if $M$ is pure of weight $m$ and $N$ is pure of weight $n$, the group of homomorphisms $ {\Hom}(M,N)$ is trivial if $ m \neq n.$
    \item  \emph{The tensor product}: there exists the tensor product $M \otimes N$ of two motives $M$ and $N$. This tensor product is compatible with the weight filtration ${\W}_*$, i.e.
\[ {\W}_n(M \otimes N) =  \sum_{i+j=n}{\W}_i(M) \otimes {\W}_j(N). \]
In terms of pure motives, if $M$ is pure of weight $m$ and $N$ is pure of weight $n$ then $ M \otimes N $ is a pure motif of weight $m+n$.
\end{enumerate}
 It is not yet clear if these properties,
which are clearly expected to be truth for motives defined over a field,
 are satisfied also by motives defined over an arbitrary scheme $S$.

If $P$, $Q$ are $S$-group schemes, we denote by $P_Q$ the fibred product $P \times_S Q$ of $P$ and $Q$ over $S$.

Let $P$, $Q$ and $G$ be commutative $S$-group schemes. A \textbf{biextension of $(P,Q)$ by $G$} is a $G_{P\times Q}$-torsor $B$, endowed with a structure of commutative extension of $Q_P$
by $G_P$ and a structure of commutative extension of $P_Q$ by $G_Q,$ which are compatible one with another (for the definition of compatible extensions see~\cite{SGA7} Expos\'e VII D\'efinition 2.1).

An \textbf{abelian $S$-scheme $A$} is an $S$-group scheme which is smooth, proper over $S$ and with connected fibres. An \textbf{$S$-torus $Y(1)$} is an $S$-group scheme which is locally isomorphic for the fpqc topology (equivalently for the \'etale topology) to an $S$-group scheme of the kind ${\GG}_m^r$ (with $r$ an integer bigger or equal to 0).
The \textbf{character group} $Y^{\vee}={\uHom}(Y(1),{\GG}_m)$ and the
\textbf{cocharacter group} $ Y={\uHom}({\GG}_m,Y(1))$
of an $S$-torus $Y(1)$
are $S$-group schemes which are \textbf{locally}
for the \'etale topology \textbf{constant group schemes defined by
 finitely generated free $\ZZ$-modules}.

 A 1-motive $M=(X,A,Y(1),G,u)$ can be viewed also as a complex $[X \stackrel{u}{\longrightarrow}G]$ of commutative $S$-group schemes concentrated in degree 0 and -1. A morphism of 1-motives is a morphism of complexes of commutative $S$-group schemes.
An \textbf{isogeny between two 1-motives}
$M_{1}=[X_{1}  \rightarrow G_{1}]$ and
$M_{2}=[X_{2} \rightarrow G_{2}]$ is a morphism of complexes $(f_{X},f_{G})$ such that
 $f_{X}:X_{1} \rightarrow X_{2}$ is injective with finite cokernel, and
$f_{G}:G_{1} \rightarrow G_{2}$ is surjective with finite kernel.
The weight filtration $W_*$ on $M=[X \rightarrow G] $ is
\begin{eqnarray}
\nonumber  {\W}_{i}(M) &=& M  ~~~~{\rm  for ~ each~} i \geq 0, \\
\nonumber  {\W}_{-1}(M) &=& [0 \longrightarrow G], \\
\nonumber  {\W}_{-2}(M) &=& [0 \longrightarrow  Y(1)], \\
\nonumber  {\W}_{j}(M) &=& 0 ~~~~{\rm  for ~ each~} j \leq -3.
\end{eqnarray}
Defining ${\rm Gr}_{i}^{{\W}}= {\W}_i / {\W}_{i+1}$, we have  ${\rm Gr}_{0}^{{\W}}(M)=
[X \rightarrow 0], {\rm Gr}_{-1}^{{\W}}(M)=
[0  \rightarrow A]$ and $ {\rm Gr}_{-2}^{{\W}}(M)=
[0  \rightarrow  Y(1)].$ Hence locally constant group schemes, abelian schemes and tori are the pure 1-motives underlying $M$ of weights 0,-1,-2 respectively. Moreover for 1-motives the weight filtration $W_*$ is defined over $\ZZ$. This means the following thing: first recall that a mixed Hodge structure $({\h}_{\ZZ}, {\W}_*, {\F}^*)$ consists of a finitely generated $\ZZ$-module ${\h}_\ZZ$, an increasing filtration $W_*$ (the weight filtration) on ${\h}_{\ZZ} \otimes \QQ$, a decreasing filtration ${\F}^*$ (the Hodge filtration) on ${\h}_{\ZZ} \otimes \CC$, and  some axioms relating these two filtrations. In the case of 1-motives defined over $\CC$, the weight filtration $W_*$ of the corresponding mixed Hodge structure (see~\cite{D1} 10.1.3) is defined on ${\h}_{\ZZ}$, and so we say that the weight filtration $W_*$ is defined over $\ZZ$.

There is a more symmetrical definition of 1-motives:
Consider the 7-tuple $(X,Y^{\vee},$\\
$A,A^*, v ,v^*,\psi)$ where
\begin{itemize}
    \item $X$ and $Y^{\vee}$ are two $S$-group schemes which are locally
for the \'etale topology constant group schemes defined by
 finitely generated free $\ZZ\,$-modules. We have to think at $X$ and at $Y^{\vee}$ as character groups of $S$-tori that we write $X^{\vee}(1)$ and $Y(1)$ and whose cocharacter groups are $X^{\vee}$ and $Y$ respectively;
    \item $A$ is an abelian $S$-scheme and $A^*$ is the dual abelian $S$-scheme of $A$ (see\\ \cite{Mu} Chapter 6 \S1);
    \item $v:X \rightarrow A$ and $v^*:Y^{\vee} \rightarrow A^*$ are two morphisms of $S$-group schemes; and
    \item $\psi$ is a trivialization of the pull-back $(v,v^*)^*{\mathcal{P}}_A$ via $(v,v^*)$ of the Poincar\' e biextension ${\mathcal{P}}_A$ of $(A,A^*)$ by ${\ZZ}(1)$.
\end{itemize}
According to Proposition~\cite{D1} (10.2.14) to have the data $(X,Y^{\vee},A,A^*, v ,v^*,\psi)$ is equivalent to have the 1-motive $M=[X {\buildrel u \over \longrightarrow }G]$, where $G$ is an extension of the abelian $S$-scheme $A$ by the $S$-torus $Y(1)$.
With these notations the Cartier dual of $M=(X,Y^{\vee},A,A^*, v ,v^*,\psi)$
is the 1-motive $M=(Y^{\vee},X,A^*,A,v^*,v, \psi \circ sym)$ where $sym: X \times Y^{\vee} \rightarrow Y^{\vee} \times X$ is the morphism which permutes the factors. \\
The pull-back $(v,v^*)^*{\mathcal{P}}_A$ by $(v,v^*)$
of the Poincar\' e biextension ${\mathcal{P}}_A$ of $(A,A^*)$ is a biextension
of $(X,Y^{\vee})$ by ${\GG}_m$.
 According~\cite{SGA3} Expos\'e X Corollary 4.5, we can suppose
that the character group $Y^{\vee}$ is $\ZZ^{\rk Y^{\vee}}$
(if necessary we localize over $S$ for the \'etale topology).
Moreover since by \cite{SGA7} Expos\'e VII (2.4.2) the category $\bBiext$ is additive in each variable, we have the equivalence of categories
\[{\bBiext}(X,Y^{\vee};{\GG}_m) \cong {\bBiext}(X,{\ZZ}; Y(1)).\]
\pn We denote by $((v,v^*)^*{\mathcal{P}}_A) \otimes Y$
 the biextension of $(X,\ZZ)$ by $Y(1)$ corresponding to the
 biextension $(v,v^*)^*{\mathcal{P}}_A$ through this equivalence of categories.
The trivialization $\psi$ of $(v,v^*)^*{\mathcal{P}}_A$ defines a trivialization
$\psi \otimes Y$ of $((v,v^*)^*{\mathcal{P}}_A) \otimes Y$, and vice versa.


\section{Biextensions of 1-motives}

\subsection{The category of biextensions of 1-motives by 1-motives}

Let $K_i=[A_i \stackrel{u_i}{\longrightarrow} B_i]$ (for $i=1,2,3$) be a complex of abelian sheaves (over any topos $\mathrm{{\mathbf{T}}}$) concentrated in degree 0 and -1.

\begin{defn}\label{defbiext}
A \textbf{biextension $({\mathcal{B}}, \Psi_1, \Psi_2,\lambda)$ of $(K_1,K_2)$ by $K_3$} consists of
\begin{enumerate}
    \item a biextension of $\mathcal{B}$ of $(B_1, B_2)$ by $B_3$;
    \item a trivialization $\Psi_1$ (resp. $\Psi_2$) of the biextension $(u_1,id_{B_2})^*{\mathcal{B}}$
    of $(A_1, B_2)$ by $B_3$ (resp. of the biextension $(id_{B_1},u_2)^*{\mathcal{B}}$ of $(B_1, A_2)$ by $B_3$) obtained as pull-back of ${\mathcal{B}}$ via $(u_1,id_{B_2}): A_1 \times B_2 \rightarrow B_1 \times B_2$ (resp. via $(id_{B_1}, u_2): B_1 \times A_2 \rightarrow B_1 \times B_2$ ). These two trivializations  have to coincide over $(A_1, A_2)$;
    \item a morphism $\lambda: A_1 \otimes A_2 \rightarrow A_3$ such that the composite $ A_1 \otimes A_2 \stackrel{\lambda}{\longrightarrow} A_3 \stackrel{u_3}{\longrightarrow} B_3$ is compatible with the restriction over $(A_1, A_2)$ of the trivializations $\Psi_1$ and $\Psi_2$.
\end{enumerate}
\end{defn}

Let $K_i=[A_i \stackrel{u_i}{\longrightarrow}B_i]$ and
$K'_i=[A'_i \stackrel{u'_i}{\longrightarrow}B'_i]$ (for $i=1,2,3$) be a complex of abelian sheaves concentrated in degree 0 and -1. Let $({\mathcal{B}},\Psi_{1}, \Psi_{2},\lambda)$ be a biextension of $(K_1,K_2)$ by $K_3$ and let $({\mathcal{B}}',\Psi'_{1}, \Psi'_{2},\lambda')$ be a biextension of $(K'_1,K'_2)$ by $K'_3$.

\begin{defn}\label{defmorphbiext}
A \textbf{morphism of biextensions}
\[(\underline{F},\underline{\Upsilon}_1,\underline{\Upsilon}_2,g_3):({\mathcal{B}},\Psi_{1}, \Psi_{2},\lambda)
\longrightarrow ({\mathcal{B}}',\Psi'_{1}, \Psi'_{2},\lambda')\]
\pn consists of
\begin{enumerate}
    \item a morphism $\underline{F}=(F,f_1,f_2,f_3):{\mathcal{B}} \rightarrow {\mathcal{B}}'$ from the biextension ${\mathcal{B}}$ to the biextension ${\mathcal{B}}'$. In particular,
\[ f_1:B_1 \longrightarrow B'_1 \qquad f_2:B_2 \longrightarrow B'_2 \qquad f_3:B_3 \longrightarrow B'_3\]
\pn are morphisms of abelian sheaves.
    \item a morphism of biextensions
\[\underline{\Upsilon}_1 = (\Upsilon_1,g_1,f_2,f_3):(u_1, id_{B_2})^* {\mathcal{B}} \longrightarrow (u'_1, id_{B'_2})^* {\mathcal{B}}'\]
\pn compatible with the morphism $\underline{F}=(F,f_1,f_2,f_3)$ and with the trivializations $\Psi_{1}$ and $\Psi'_{1}$, and  a morphism of biextensions
\[\underline{\Upsilon}_2 = (\Upsilon_2,f_1,g_2,f_3):(id_{B_1},u_2)^* {\mathcal{B}} \longrightarrow (id_{B'_1},u'_2)^* {\mathcal{B}}'\]
\pn compatible with the morphism $\underline{F}=(F,f_1,f_2,f_3)$ and with the trivializations $\Psi_{2}$ and $\Psi'_{2}$. In particular,
\[g_1: A_1 \longrightarrow A'_1 \qquad g_2: A_2 \longrightarrow A'_2  \]
\pn are morphisms of abelian sheaves. By pull-back, the two morphisms \\ $\underline{\Upsilon}_1=(\Upsilon_1,g_1,f_2,f_3)$ and $\underline{\Upsilon}_2=(\Upsilon_2,f_1,g_2,f_3)$ define a morphism of biextensions
 $\underline{\Upsilon}=(\Upsilon,g_1,g_2,f_3): (u_1,u_2)^* {\mathcal{B}} \rightarrow (u'_1, u'_2)^* {\mathcal{B}}'$ compatible with the morphism $\underline{F}=(F,f_1,f_2,f_3)$ and with the trivializations $\Psi$ and $\Psi'$.
    \item a morphism $g_3:A_3 \rightarrow A'_3$ of abelian sheaves compatible with $u_3$ and $u'_3$ (i.e. $u'_3 \circ g_3 =f_3 \circ u_3$) and such that
\[ \lambda' \circ (g_1 \times g_2)= g_3 \circ \lambda\]
\end{enumerate}
\end{defn}

\begin{rem}\label{remmorphbiext}
 The morphisms $g_3$ and $f_3$ define a morphism from $K_3$ to $K'_3$. The morphisms $g_1$ and $f_1$ (resp. $g_2$ and $ f_2$) define morphisms from $K_1$ to $K'_1$ (resp. from $K_2$ to $K'_2$).
\end{rem}

We denote by ${\bBiext}(K_1,K_2;K_3)$ the category of biextensions
of $(K_1,K_2)$ by $K_3$. The Baer sum of extensions defines a group law for the objects of the category ${\bBiext}(K_1,K_2;K_3)$, which is therefore a Picard category (see~\cite{SGA7} Expos\'e VII 2.4, 2.5 and 2.6).
Let ${\Biext}^0(K_1,K_2;K_3)$ be the group of automorphisms of any biextension of $(K_1,K_2)$ by $K_3$, and let\\
 ${\Biext}^1(K_1,K_2;K_3)$ be the group of isomorphism classes of biextensions of $(K_1,K_2)$ by $K_3$.

\begin{rem}\label{homolo} In the paper~\cite{B4} in preparation, we are proving the following homological interpretation of the groups ${\Biext}^i(K_1,K_2;K_3)$ for $i=0,1$:
\[ {\Biext}^i(K_1,K_2;K_3) \cong {\Ext}^i(K_1 {\buildrel {\scriptscriptstyle \LL} \over \otimes} K_2,K_3) \qquad (i=0,1). \]
This homological interpretation generalizes the one obtained by Grothendick in\\ \cite{SGA7} Expos\'e VII 3.6.5 for biextensions of abelian sheaves.
\end{rem}

Now we generalize to complex of abelian sheaves concentrated in degree 0 and -1, the definitions of symmetric biextensions and of skew-symmetric biextensions of abelian sheaves introduced by L. Breen
in~\cite{Be1} 1.4 and~\cite{Be2} 1.1 respectively. Let $K$ and $K'$ be complexes of abelian sheaves concentrated in degree 0 and -1. Denote by $sym:K \times K \rightarrow K \times K$ the morphism which permutes the factors and by
 $d:K \rightarrow K \times K$ the diagonal morphism.

\begin{defn} \label{biextsym}
A \textbf{symmetric biextension $({\mathcal{B}},\xi_{\mathcal{B}})$
of $(K,K)$ by $K'$} consists of a biextension ${\mathcal{B}}=({\mathcal{B}},\Psi_{1}, \Psi_{2},\lambda)$ of $(K,K)$ by $K'$
and a morphism of biextensions
$ \xi_{\mathcal{B}}: {sym}^*{\mathcal{B}} \rightarrow {\mathcal{B}},$ where ${sym}^*{\mathcal{B}}$ is the pull-back of ${\mathcal{B}}$ via the morphism $sym$ which permutes the factors, such that the restriction $d^* \xi_{\mathcal{B}}$ of $ \xi_{\mathcal{B}}$ by the diagonal morphism $d$ coincides with the isomorphism
$$\nu_{\mathcal{B}} : d^*{sym}^* {\mathcal{B}} \longrightarrow d^*{\mathcal{B}}$$
arising from the identity $sym\circ  d =d.$
\end{defn}

The morphism $\xi_{\mathcal{B}}$ is involute, i.e. the composite
 $\xi_{\mathcal{B}} \circ {sym}^*\xi_{\mathcal{B}}:
{sym}^*{sym}^* {\mathcal{B}} \rightarrow {sym}^*{\mathcal{B}} \rightarrow {\mathcal{B}}$  is the identity of ${\mathcal{B}}$ (cf. [Br83] 1.7).

\begin{defn}
\textbf{The symmetrized biextension} of a biextension ${\mathcal{B}}=({\mathcal{B}},\Psi_{1},$\\
$ \Psi_{2},\lambda)$ of $(K,K)$ by $K'$ is the symmetric biextension $({\mathcal{B}} \wedge {sym}^* {\mathcal{B}}, \xi_{{\mathcal{B}} \wedge {sym}^* {\mathcal{B}}})$, where the morphism   $ \xi_{{\mathcal{B}} \wedge {sym}^* {\mathcal{B}}}$ is given canonically by the composite
$$\xi_{{\mathcal{B}} \wedge {sym}^* {\mathcal{B}}}: {sym}^* {\mathcal{B}} \wedge {sym}^* {sym}^* {\mathcal{B}} \longrightarrow {sym}^* {\mathcal{B}} \wedge  {\mathcal{B}}
{\buildrel \tau \over \longrightarrow} {\mathcal{B}} \wedge {sym}^* {\mathcal{B}}$$
where the first arrow comes from the equality $ sym \circ   sym =  id$ and the second one is the morphism
$\tau:{sym}^* {\mathcal{B}} \wedge {\mathcal{B}} \rightarrow  {\mathcal{B}} \wedge {sym}^* {\mathcal{B}}$ which permutes the factors of the contracted product.
\end{defn}

\begin{defn}\label{antisym}
A \textbf{skew-symmetric biextension $({\mathcal{B}},\varphi_{\mathcal{B}})$
of $(K,K)$ by $K'$} consists of a biextension ${\mathcal{B}}=({\mathcal{B}},\Psi_{1}, \Psi_{2},\lambda)$ of $(K,K)$ by $K'$
and a trivialization $\varphi_{\mathcal{B}}$ of the biextension $({\mathcal{B}} \wedge {sym}^* {\mathcal{B}}, \xi_{{\mathcal{B}} \wedge {sym}^* {\mathcal{B}}})$ which is compatible with the symmetric structure
of $({\mathcal{B}} \wedge {sym}^* {\mathcal{B}}, \xi_{{\mathcal{B}} \wedge {sym}^* {\mathcal{B}}})$.
\end{defn}

Since we can view 1-motives as
complexes of commutative $S$-group schemes concentrated in degree 0 and -1, all the definitions of this section apply to 1-motives.

\subsection{A simpler description}
  Using the symmetrical description of 1-motives recalled in the Notation, we can give a simpler description of the definition~\ref{biext1-mot} of biextension of 1-motives by 1-motives.

\begin{prop}\label{defbiext2}
Let $M_i=(X_i,Y_i^{\vee},A_i,A_i^*, v_i ,v_i^*,\psi_i)$ (for $i=1,2,3$) be a 1-motive over $S$. To have a biextension $({\mathcal{B}},\Psi_{1}, \Psi_{2},\lambda)$ of $(M_1,M_2)$ by $M_3$ is equivalent to have a 4-uplet $(B,\Phi_1, \Phi_2,\Lambda)$ where
\begin{enumerate}
    \item a biextension of $B$ of $(A_1, A_2)$ by $Y_3(1);$
    \item a trivialization $\Phi_1 $ (resp. $\Phi_2$) of the biextension
$(v_1, id_{A_2})^*B$ of $(X_1, A_2)$ by $Y_3(1)$  (resp. of the biextension
 $(id_{A_1},v_2)^*B$ of $(A_1, X_2)$ by $Y_3(1)$) obtained as pull-back of  $B$ via $(v_1,id_{A_2})$ (resp. via $(id_{A_1},v_2)$).
These two trivializations $\Phi_1 $ and $\Phi_2 $ have to coincide over $(X_1, X_2)$, i.e.
\[(v_1,id_{X_2})^*\Phi_2=\Phi=(id_{X_1},v_2)^*\Phi_1 \]
\pn with $\Phi$ a trivialization of the biextension $(v_1, v_2)^*{B}$
 of $(X_1, X_2)$ by $Y_3(1)$ obtained as pull-back of the biextension
$B$ via $(v_1,v_2)$;
    \item a morphism $\Lambda: (v_1, v_2)^*B \rightarrow ((v_3,v_3^*)^*{\mathcal{P}}_{A_3}) \otimes Y_3 $ of biextensions, with $\Lambda_{| Y_3(1)}$ equal to the identity and $\Lambda_{| X_1 \times X_2}$ bilinear, such that the following diagram is commutative
\begin{equation}\label{condbiext}
\begin{array}{ccc}
  Y_3(1) &= &  Y_3(1)\\
        \vert &     & \vert\\
(v_1, v_2)^*B &  \longrightarrow &((v_3,v_3^*)^*{\mathcal{P}}_{A_3}) \otimes Y_3 \\
{\scriptstyle \Phi} \uparrow \downarrow~~ & &~~~~~~~~~\downarrow \uparrow {\scriptstyle \psi_3 \otimes Y_3}   \\
X_1 \times X_2 & \longrightarrow & X_3 \times \ZZ. \\
\end{array}
\end{equation}
\end{enumerate}
\end{prop}

\begin{proof}
According to~\cite{B3} Theorem 2.5.2 and remark 2.5.3, to have the biextension $B$ of $(A_1,A_2)$ by $Y_3(1)$ is equivalent to have the  biextension
${\mathcal{B}}=\iota_{3\,*}(\pi_1,\pi_2)^*B$ of $(G_1, G_2)$ by $G_3$, where for $i=1,2,3,$ $\pi_i:G_i \rightarrow A_i$
is the projection of $G_i$ over $A_i$ and $\iota_i:Y_i(1) \rightarrow
G_i$ is the inclusion of $Y_i(1)$ over $G_i.$
The trivializations $(\Phi_1 ,\Phi_2 )$ and $(\Psi_1 ,\Psi_2 )$ determine each others.
To have the morphism of $S$-group schemes $\lambda:X_1 \times X_2 \rightarrow X_3$ is equivalent to have the morphism
of biextensions
$\Lambda: (v_1, v_2)^*B \rightarrow ((v_3,v_3^*)^*{\mathcal{P}}_{A_3}) \otimes Y_3 $  with $\Lambda_{| Y_3(1)}$ equal to the identity. In particular, through this last equivalence $\lambda$ corresponds to
$\Lambda_{| X_1 \times X_2}$ and
to require that
$u_3 \circ \lambda: X_1 \otimes X_2 \rightarrow G_3$ is compatible with the trivialization  $\Psi$ of $(u_1,u_2)^* {\mathcal{B}}$ corresponds to require the commutativity of the diagram (\ref{condbiext}) with the vertical arrows going up.
\end{proof}


\section{Some tensor products}

\subsection{The tensor product with a motive of weight zero}

 Let $Z$ be an $S$-group scheme which is locally
for the \'etale topology a constant group scheme defined by a
 finitely generated free $\ZZ$-module, i.e.
 there exist an \'etale surjective morphism $S' \rightarrow S$ such that
 $Z_{S'}$ is the constant $S'$-group scheme ${\ZZ}^z$ with $z$ an integer bigger or equal to 0. The tensor product of abelian sheaves in the big \'etale site furnishes the tensor product of $Z$ with the pure motives underlying 1-motives. In this section we discuss the representability by group schemes of such tensor products.

\subsubsection{The tensor product of two motives of weight zero:}
Let $X$ be an $S$-group scheme which is locally
for the \'etale topology a constant group scheme defined by a
 finitely generated free $\ZZ$-module, ${\ZZ}^x$ with $x$ an integer bigger or equal to 0. The tensor product
$$X \otimes Z$$
 is the $S$-group scheme which is locally
for the \'etale topology a constant group scheme defined by a
 finitely generated free $\ZZ$-module of rank $x \cdot z$, such that
  there exist an \'etale surjective morphism $S' \rightarrow S$ for which
 the $S'$-group scheme $(X \otimes Z)_{S'}$ is isomorphic to the fibred product of $z$-copies of the $S'$-group scheme $X_{S'}.$
In fact, let $g: S' \rightarrow S$ be the \'etale surjective morphism such that $Z_{S'}$ is the constant $S'$-group scheme ${\ZZ}^z$. Over $S'$ we define the tensor product $X_{S'} \otimes Z_{S'}$ as the fibred product of  $z$-copies of the ${S'}$-scheme $X_{S'}$:
$$X_{S'} \otimes Z_{S'} := X_{S'} \times_{S'} \ldots \times_{S'} X_{S'}$$
The $S'$-scheme $X_{S'} \otimes Z_{S'}$ is again an $S'$-group scheme which is locally for the \'etale topology a constant group scheme defined by a
 finitely generated free $\ZZ$-module of rank $x \cdot z$, and so in particular it is locally of finite presentation, separated and locally quasi-finite over $S'$. By~\cite{SGA3} Expos\'e X Lemma 5.4 the morphism $g$  is a morphism of effective descent for the fibred category of locally of finite presentation, separated and locally quasi-finite schemes, i.e. there exists an $S$-scheme $X \otimes Z$ and an $S'$-isomorphism $(X \otimes Z)_{S'} \cong X_{S'} \otimes Z_{S'}$
which is compatible with the descent data. By construction $X \otimes Z$ is again a group scheme which is locally for the \'etale topology a constant group scheme defined by a
 finitely generated free $\ZZ$-module of rank $x \cdot z$.

\subsubsection{The tensor product of a torus with a motive of weight 0:} Let $Y(1)$ be an $S$-torus.
The tensor product
$$Y(1) \otimes Z$$
 is the $S$-torus such that
  there exist an \'etale surjective morphism $S' \rightarrow S$ for which
 the $S'$-torus $(Y(1) \otimes Z)_{S'}$ is isomorphic to the fibred product of $z$-copies of the $S'$-torus $Y(1)_{S'}.$ In fact, let $g: S' \rightarrow S$ be the \'etale surjective morphism such that $Z_{S'}$ is the constant $S'$-group scheme ${\ZZ}^z$. Over $S'$ we define the tensor product $Y(1)_{S'} \otimes Z_{S'}$ as the fibred product of  $z$-copies of the ${S'}$-torus $Y(1)_{S'}$:
$$Y(1)_{S'} \otimes Z_{S'} := Y(1)_{S'} \times_{S'} \ldots \times_{S'} Y(1)_{S'}$$
The $S'$-scheme $Y(1)_{S'} \otimes Z_{S'}$ is again an $S'$-torus, and so in particular it is affine $S'$. By~\cite{SGA1} Expos\'e VIII Theorem 2.1 the morphism $g$ is a morphism of effective descent for the fibred category of affine schemes, i.e. there exists an $S$-scheme $Y(1) \otimes Z$ and an $S'$-isomorphism $(Y(1) \otimes Z)_{S'} \cong Y(1)_{S'} \otimes Z_{S'}$
which is compatible with the descent data. By construction $Y(1) \otimes Z$ is again a torus.

\begin{rem} If the cocharacter group $Y$ of the torus $Y(1)$ has rank $y$,
the cocha\-ra\-cter group of the torus $Y(1) \otimes Z$ is the motif of weight zero $Y \otimes Z$ of rank $y \cdot z$.
\end{rem}

\subsubsection{The tensor product of an abelian scheme with a motive of weight 0:} Let $A$ be an abelian $S$-scheme.
The tensor product
$$A \otimes Z$$ is the abelian $S$-scheme  such that
  there exist an \'etale surjective morphism $S' \rightarrow S$ for which
 the abelian $S'$-scheme $(A \otimes Z)_{S'}$ is isomorphic to the fibred product of $z$-copies of the abelian $S'$-scheme $A_{S'}.$
In fact, let $g: S' \rightarrow S$ be the \'etale surjective morphism such that $Z_{S'}$ is the constant $S'$-group scheme ${\ZZ}^z$. Over $S'$ we define the tensor product $A_{S'} \otimes Z_{S'}$ as the fibred product of  $z$-copies of the abelian ${S'}$-scheme $A_{S'}$:
$$A_{S'} \otimes Z_{S'} := A_{S'} \times_{S'} \ldots \times_{S'} A_{S'}$$
The $S'$-scheme $A_{S'} \otimes Z_{S'}$ is again an abelian $S'$-scheme, and in particular it is an algebraic space over $S'$. By~\cite{LM-B} Corollary 10.4.2 the morphism $g$ is a morphism of effective descent for the fibred category of algebraic spaces, i.e. there exists an algebraic $S$-space $A \otimes Z$ and an $S'$-isomorphism $(A \otimes Z)_{S'} \cong A_{S'} \otimes Z_{S'}$ which is compatible with the descent data.
The local properties, as smoothness, and the properties which are stable by base change, as properness and geometrically connected fibres, are carried over from $A$ to $A \otimes Z$. Therefore
the algebraic $S$-space $A \otimes Z$ is a group object, smooth, proper over $S$ and with connected fibres and so according to~\cite{F} Theorem 1.9, $A \otimes Z$ is an abelian $S$-scheme.

\subsubsection{The tensor product of an extension of an abelian scheme by a torus with a motive of weight 0:} Let $G$ be an extension of an abelian $S$-scheme $A$ by an $S$-torus $Y(1)$. The tensor product
$$G \otimes Z$$
 is the $S$-group scheme which is extension of the abelian $S$-scheme $A \otimes Z$ by the $S$-torus $Y(1) \otimes Z$, such that
  there exist an \'etale surjective morphism $S' \rightarrow S$ for which
 the $S'$-group scheme $(G \otimes Z)_{S'}$ is isomorphic to the fibred product of $z$-copies of the $S'$-group scheme $G_{S'}.$
In fact, let $g: S' \rightarrow S$ be the \'etale surjective morphism such that $Z_{S'}$ is the constant $S'$-group scheme ${\ZZ}^z$. Over $S'$ we define the tensor product $G_{S'} \otimes Z_{S'}$ as the fibred product of  $z$-copies of the ${S'}$-scheme $G_{S'}$:
$$G_{S'} \otimes Z_{S'} := G_{S'} \times_{S'} \ldots \times_{S'} G_{S'}$$
The $S'$-scheme $G_{S'} \otimes Z_{S'}$ is an extension of the abelian $S'$-scheme $A_{S'} \otimes Z_{S'}$ by an $S'$-torus $Y(1)_{S'} \otimes Z_{S'}.$ In particular $G_{S'} \otimes Z_{S'}$ is an algebraic space over $S'$. By~\cite{LM-B} Corollary 10.4.2 the morphism $g$ is a morphism of effective descent for the fibred category of algebraic spaces, i.e. there exists an algebraic $S$-space $G \otimes Z$ which is
an extension the abelian $S$-scheme $A \otimes Z$ by the $S$-torus $Y(1) \otimes Z$, and an $S'$-isomorphism $(G \otimes Z)_{S'} \cong G_{S'} \otimes Z_{S'}$ which is compatible with the descent data.
The smoothness of the torus $Y(1) \otimes Z$ implies that the fppf $Y(1) \otimes Z$-torsor $G \otimes Z$ is in fact an \'etale torsor (see~\cite{G3} Theorem 11.7).
Since the torus $Y(1) \otimes Z$ is affine over $S$ and since affiness is stable under base extensions, the $Y(1) \otimes Z$-torsor
$G \otimes Z$ is affine over $A \otimes Z$.
By the theory of effective descent for the fibred category of affine schemes (see~\cite{SGA1} Expos\'e VIII Theorem 2.1) the algebraic $S$-space $G \otimes Z$ is in fact an $S$-scheme.

Using the above constructions we can now define the tensor product of a 1-motive with a motive of weight 0:

\begin{defn} Let $M=[X \stackrel{u}{\rightarrow} G]$ be a 1-motive.
The tensor product $M \otimes Z$ is the 1-motive
\[[X \otimes Z \, \stackrel{u \otimes Z}{\longrightarrow} \, G \otimes Z]. \]
\end{defn}

We can conclude that roughly speaking ``to tensor a motive by a motive of weight 0'' means to take a certain number of copies of this motive.

\subsection{The 1-motive underlying $M_1 \otimes M_2 / {\W}_{-3}(M_1 \otimes M_2)$} Let $M_i=(X_i,Y^{\vee}_i,A_i,$\\
$ A^*_i,v_i ,v^*_i,\psi_i)$ be a 1-motive (for $i=1,2$) defined over $S$. The 1-motive underlying the motive
$M_1 \otimes M_2 / {\W}_{-3}(M_1 \otimes M_2)$ is the 1-motive
$\mathbb{M}=(\mathbb{X},\mathbb{Y}^{\vee},\mathbb{A},\mathbb{A}^*,\mathbb{V} ,\mathbb{V}^*,\Psi)$ where,
\begin{itemize}
    \item $\mathbb{X}$ is the $S$-group scheme $X_1 \otimes X_2$,
    \item $\mathbb{Y}^{\vee}$ is the $S$-group scheme $X_1^{\vee} \otimes Y_2^{\vee} ~ + ~ Y_1^{\vee} \otimes X_2^{\vee} ~ + ~ {\mathcal{B}iext}^1(A_1, A_2; {\ZZ}(1)) $ where ${\mathcal{B}iext}^1(A_1, A_2; {\ZZ}(1))$ is the $S$-group scheme  which is locally
for the \'etale topology a constant group scheme defined by the group
of isomorphism classes of biextensions of $(A_1, A_2)$ by ${\ZZ}(1)$ (remark that ${\mathcal{B}iext}^1(A_1, A_2; {\ZZ}(1))$ is a group scheme because ${\Biext}^1(A_1, A_2; {\ZZ}(1))$ is a group),
    \item $\mathbb{A}$ is the abelian $S$-scheme $X_1 \otimes A_2 ~ + ~ A_1 \otimes X_2$,
    \item the morphisms $\mathbb{V}$ and $\mathbb{V}^*$ and the trivialization $\Psi$ are defined by the formula~(\ref{V}), (\ref{Vdual-1}),  (\ref{Vdual-2}), (\ref{Vdual-3}), (\ref{Psi-1}), (\ref{Psi-2}), (\ref{Psi-3}).
\end{itemize}

\begin{proof}
The only non trivial components of the motive $M_1 \otimes M_2 / {\W}_{-3}(M_1 \otimes M_2)$ are the pure motives
\begin{eqnarray}
\nonumber  {\Gr}_{0}^W(M_1 \otimes M_2 / {\W}_{-3}(M_1 \otimes M_2)) &=& X_1 \otimes X_2, \\
\nonumber  {\Gr}_{-1}^W(M_1 \otimes M_2 / {\W}_{-3}(M_1 \otimes M_2)) &=& X_1 \otimes A_2 ~ +~ A_1 \otimes X_2, \\
\nonumber  {\Gr}_{-2}^W(M_1 \otimes M_2 / {\W}_{-3}(M_1 \otimes M_2)) &=& X_1 \otimes Y_2(1) ~ + ~ Y_1(1) \otimes X_2 ~ + ~A_1 \otimes A_2.
\end{eqnarray}
of weight 0,-1 and -2 respectively. Until now we don't have defined $A_1 \otimes A_2$ but in fact, in what follows, we will need only the morphisms from $A_1 \otimes A_2$ to the torus ${\ZZ}(1)$ which are defined in Definition~\ref{defmor} (see also remark~\ref{rem}). Hence
\begin{itemize}
    \item $\mathbb{X}$ is the $S$-group scheme $X_1 \otimes X_2$ which is locally
for the \'etale topology a constant group scheme defined by a finitely generated free $\ZZ$-module of rank $r_1 \cdot r_2$, where $r_1$ (resp. $r_2$) is the rank of the finitely generated free $\ZZ$-module defining $X_1$ (resp.$X_2$),
    \item $\mathbb{A}$ is the abelian $S$-scheme $X_1 \otimes A_2 ~ + ~ A_1 \otimes X_2$,
    \item $\mathbb{Y}^{\vee}$ is ${\uHom}(X_1 \otimes Y_2(1) ~ + ~ Y_1(1) \otimes X_2 ~ + ~A_1 \otimes A_2, {\ZZ}(1))$.
     \end{itemize}
We have the equality
 $$\mathbb{Y}^{\vee} =  X_1^{\vee} \otimes Y_2^{\vee} ~ + ~ Y_1^{\vee} \otimes X_2^{\vee} ~ + ~ {\uHom}(A_1 \otimes A_2, {\ZZ}(1)). $$
As we will see in definition~\ref{defmor}, the bilinear morphisms from
 $A_1 \times A_2$ to the torus ${\ZZ}(1)$ are the isomorphism classes of biextensions of $(A_1, A_2)$ by ${\ZZ}(1)$. Therefore
$$\mathbb{Y}^{\vee} =  X_1^{\vee} \otimes Y_2^{\vee} ~ + ~ Y_1^{\vee} \otimes X_2^{\vee} ~ + ~ {\mathcal{B}iext}^1(A_1, A_2; {\ZZ}(1)) $$
where ${\mathcal{B}iext}^1(A_1, A_2; {\ZZ}(1))$ is the $S$-group scheme  which is locally
for the \'etale topology a constant group scheme defined by the group
${\Biext}^1(A_1, A_2; {\ZZ}(1))$ of isomorphism classes of biextensions of $(A_1, A_2)$ by ${\ZZ}(1)$.\\
We define the morphism $\mathbb{V}$ using the morphisms $v_1:X_1 \rightarrow A_1$ and $v_2:X_2 \rightarrow A_2$. In fact,
\begin{equation}\label{V}
\mathbb{V}=(v_1^{\otimes X_2},v_2^{\otimes X_1}): X_1 \otimes X_2 \longrightarrow  X_1 \otimes A_2 ~ + ~ A_1 \otimes X_2
\end{equation}
where $v_1^{\otimes X_2}: X_1 \otimes X_2 \rightarrow  A_1 \otimes X_2 $
  and $v_2^{\otimes X_1}: X_1 \otimes X_2 \rightarrow  X_1 \otimes A_2. $
Before to define the morphism $\mathbb{V}^*$, observe that the Cartier dual of the abelian $S$-scheme $X_1 \otimes A_2 ~ + ~ A_1 \otimes X_2$
is the abelian $S$-scheme $X_1^{\vee} \otimes A_2^* ~ + ~ A_1^* \otimes X_2^{\vee}$. Since $\mathbb{Y}^{\vee}$ decomposes in three terms, we can define $\mathbb{V}$ separately over each terms. For the first two terms $X_1^{\vee} \otimes Y_2^{\vee} ~ + ~ Y_1^{\vee} \otimes X_2^{\vee} $, we use the morphisms $v_1^*:Y_1^{\vee} \rightarrow A_1^*$ and $v_2^*: Y_2^{\vee} \rightarrow A_2^*$. In fact,
\begin{equation}\label{Vdual-1}
    \mathbb{V}^*_{1} =(v_1^*)^{ \otimes X_2^{\vee}}: Y_1^{\vee} \otimes X_2^{\vee}   \longrightarrow  A_1^* \otimes X_2^{\vee}
\end{equation}
\begin{equation}\label{Vdual-2}
\mathbb{V}^*_{2} = (v_2^*)^{ \otimes X_1^{\vee}}:
 X_1^{\vee} \otimes Y_2^{\vee} \longrightarrow  X_1^{\vee} \otimes A_2^*
 \end{equation}
In order to define $\mathbb{V}^*$ over the term ${\mathcal{B}iext}^1(A_1, A_2; {\ZZ}(1))$ we use the well-known canonical isomorphisms
\[{\Hom}(A_1,A_2^*) \cong {\Biext}^1(A_1,A_2;{\ZZ}(1)) \cong
{\Hom}(A_2,A_1^*)\]
(see~\cite{SGA7} Expos\'e VIII 3.2) and the isomorphisms $X_1^{\vee} \otimes A_2^* \cong {\uHom}(X_1,A_2^*)$ and $A_1^* \otimes X_2^{\vee} \cong {\uHom}(X_2,A_1^*)$:
\begin{eqnarray}\label{Vdual-3}
  \mathbb{V}^*_3:{\mathcal{B}iext}^1(A_1, A_2; {\ZZ}(1))  &\longrightarrow & {\uHom}(X_1,A_2^*) +  {\uHom}(X_2,A_1^*) \\
 \nonumber b &\longmapsto & (b \circ v_1, b \circ v_2).
\end{eqnarray}
We set $\mathbb{V}^*=\mathbb{V}^*_1 + \mathbb{V}^*_2 + \mathbb{V}^*_3.$
It remains to define the trivialization $\Psi$
 of the pull-back via $(\mathbb{V},\mathbb{V}^*)$ of the Poincar\' e biextension of $(X_1 \otimes A_2 ~ +~ A_1 \otimes X_2,X_1^{\vee} \otimes A_2^* ~ + ~ A_1^* \otimes X_2^{\vee})$ by ${\ZZ}(1)$.
Since $\mathbb{Y}^{\vee}$ decomposes in three terms, we can define this trivialization separately over each terms. For the first two terms $X_1^{\vee} \otimes Y_2^{\vee} ~ + ~ Y_1^{\vee} \otimes X_2^{\vee} $, we use the trivializations $\psi_1: X_1 \times Y_1^{\vee} \rightarrow {\ZZ}(1)$ and $\psi_2: X_2 \times Y_2^{\vee} \rightarrow {\ZZ}(1)$. In fact, the trivializations
$$ (\psi_1)^{ \otimes X_2}: (X_1 ~ \times ~   Y_1^{\vee}) \otimes X_2  \longrightarrow {\ZZ}(1) \otimes X_2$$
$$ (\psi_2)^{ \otimes X_1}: X_1 \otimes (X_2  ~ \times ~ Y_2^{\vee})
  \longrightarrow X_1 \otimes {\ZZ}(1) $$
furnish
\begin{equation}\label{Psi-1}
    \Psi_{1} =(\psi_1)^{ \otimes X_2} \otimes X_2^{\vee}:(X_1 \otimes X_2)  ~ \times ~   (Y_1^{\vee} \otimes X_2^{\vee} )  \longrightarrow {\ZZ}(1)
\end{equation}
\begin{equation}\label{Psi-2}
\Psi_{2} = (\psi_2)^{ \otimes X_1}\otimes X_1^{\vee}: (X_1 \otimes X_2)  ~ \times ~ (X_1^{\vee} \otimes Y_2^{\vee}) \longrightarrow {\ZZ}(1)
 \end{equation}
In order to define $\Psi$ over the term ${\mathcal{B}iext}^1(A_1, A_2; {\ZZ}(1))$ we use the definition~\ref{defmor} according to which
the group ${\Biext}^1(A_1, A_2; {\ZZ}(1))$ of isomorphism classes of biextensions of $(A_1, A_2)$ by ${\ZZ}(1)$ is the group ${\Hom}(A_1 \otimes A_2, {\ZZ}(1))$ of bilinear morphisms from $A_1 \times A_2$ to ${\ZZ}(1)$:
\begin{eqnarray}\label{Psi-3}
  \Psi_3: (X_1 \otimes X_2)  ~ \times ~ {\mathcal{B}iext}^1(A_1, A_2; {\ZZ}(1))  &\longrightarrow & {\ZZ}(1) \\
 \nonumber (x_1 \otimes x_2, b) &\longmapsto & b( v_1(x_1)\otimes v_2(x_2) ).
\end{eqnarray}
We set $\Psi=\Psi_1 + \Psi_2 + \Psi_3.$
\end{proof}

\begin{rem}\label{rem} Let $S$ be the spectrum of a field of characteristic 0 embeddable in $\CC$. The Hodge realization $T_{\h}(A_1 \otimes A_2)$ of the motive
$A_1 \otimes A_2$ is a pure Hodge structure of type $\{(-2,0),(-1,-1),(0,-2)\}$.
Since the Hodge realization of a 1-motive is of type $\{(0,0),(-1,0),(0,-1),(-1,-1)\},$ the only component of
$T_{\h}(A_1 \otimes A_2)$ which is involved in the Hodge realization of the 1-motive underlying $M_1 \otimes M_2 / {\W}_{-3}(M_1 \otimes M_2),$ is the component $(T_{\h}(A_1 \otimes A_2))^{-1,-1}$ of type $(-1,-1)$. By definition~\ref{defmor} the group ${\Hom}(A_1 \otimes A_2, {\ZZ}(1))$ is the group
${\Biext}^1(A_1, A_2; {\ZZ}(1))$ of isomorphism classes of biextensions of $(A_1,A_2)$ by ${\ZZ}(1)$, and so the component $\big(T_{\h}(A_1 \otimes A_2)\big)^{-1,-1}$ is the Hodge realization of the torus
$\big({\Biext}^1(A_1, A_2; {\ZZ}(1))\big)^{\vee}(1)$ whose character group is ${\Biext}^1(A_1, A_2; {\ZZ}(1))$:
\[\big(T_{\h}(A_1 \otimes A_2)\big)^{-1,-1} =T_{\h}\big(({\Biext}^1(A_1, A_2; {\ZZ}(1)))^{\vee}(1)\big).\]
\end{rem}

\section{Morphisms from a finite tensor product of 1-motives to a 1-motive}

\subsection{Multilinear morphisms between 1-motives}

\begin{defn}\label{defmor}
Let $M_1, M_2$ and $M_3$ be three 1-motives defined over $S$.
\textbf{A morphism from the tensor product of $M_1$ and $M_2$
to $M_3$} is an isomorphism class of biextensions of $(M_1, M_2)$
by $M_3$. Moreover to the three 1-motives $M_1, M_2$ and $M_3$ we associate
a group ${\Hom}(M_1, M_2;M_3)$ defined in the following way:
\[ {\Hom}(M_1, M_2;M_3) := {\Biext}^1(M_1,M_2;M_3), \]
i.e. ${\Hom}(M_1, M_2;M_3)$ is \emph{the group of bilinear morphisms from
$ M_1 \times M_2$ to $M_3$.}
\end{defn}

The structure of commutative group of ${\Hom}(M_1,M_2;M_3)$ is described in~\cite{SGA7} Expos\'e VII 2.5.

Let $M_i$ be a 1-motive $[X_i \rightarrow 0]$ of weight 0 (for $i=1,2,3$). According to our definition of biextension of 1-motives by 1-motives,
 we have the equality
 \[{\Biext}^1([X_1 \rightarrow 0],[X_2 \rightarrow 0];[X_3 \rightarrow 0]) = {\Hom}(X_1 \otimes X_2,X_3),\]
i.e. biextensions of $([X_1 \rightarrow 0],[X_2 \rightarrow 0])$ by $[X_3 \rightarrow 0]$ are just bilinear morphisms of $S$-group schemes from $X_1 \times X_2$ to $X_3$. As expected, for motives of weight 0 we have therefore
\[ {\Hom}([X_1 \rightarrow 0],[X_2 \rightarrow 0];[X_3 \rightarrow 0]) = {\Hom}(X_1 \otimes X_2,X_3). \]

Definition \ref{defmorphbiext} of morphisms of biextensions of 1-motives
by 1-motives allows us to define a morphism between the bilinear morphisms corresponding to such biextensions. More precisely,
let $M_i$ and $M'_i$ (for $i=1,2,3$) be 1-motives over $S$.
If we denote by $b$ the morphism
$ M_1 \otimes M_2 \rightarrow M_3$ corresponding to the biextension
$({{\mathcal{B}}},\Psi_{1}, \Psi_{2},\lambda)$ of $(M_1,M_2)$ by $M_3$
and by $b'$ the morphism
$ M'_1 \otimes M'_2 \rightarrow  M'_3$ corresponding to the biextension
$({{\mathcal{B}}}',\Psi'_{1}, \Psi'_{2},\lambda')$ of $(M'_1,M'_2)$ by $M'_3$,
a morphism
$(\underline{F},\underline{\Upsilon}_1,\underline{\Upsilon}_2,g_3):
({{\mathcal{B}}},\Psi_{1}, \Psi_{2},\lambda)
\rightarrow ({{\mathcal{B}}}',\Psi'_{1}, \Psi'_{2},\lambda')$
of biextensions defines the vertical arrows of the following diagram of morphisms
\[\begin{array}{ccc}
  M_1 \otimes M_2  & \stackrel{b}{\longrightarrow} &M_3  \\
 \downarrow  &  &  \downarrow \\
  M'_1 \otimes M'_2 & \stackrel{b'}{\longrightarrow} & M'_3. \\
\end{array}\]
\pn It is clear now
why from the data $(\underline{F},\underline{\Upsilon}_1,\underline{\Upsilon}_2,g_3)$ we get a
morphism from $M_3$ to $M'_3$ as
remarked in \ref{remmorphbiext}. Moreover
since $ M_1 \otimes [{\ZZ} \rightarrow 0]$,
$ M'_1 \otimes [{\ZZ} \rightarrow 0]$,
$ [{\ZZ} \rightarrow 0] \otimes  M_2$ and
$[{\ZZ} \rightarrow 0] \otimes M'_2$, are sub-1-motives
of the motives $M_1 \otimes M_2$ and $M'_1 \otimes M'_2$,
from the data
$(\underline{F},\underline{\Upsilon}_1,\underline{\Upsilon}_2,g_3)$ we get also morphisms
 from $M_1$ to $M'_1$ and from $M_2$ to $M'_2$ (see \ref{remmorphbiext}).

In~\cite{B3} Theorem 2.5.2 we proved that if $M_i=[X_i \stackrel{u_i}{\longrightarrow}G_i]$ (for $i=1,2,3$) is a 1-motive defined over $S$, ${\Biext}^1(G_1,G_2;G_3) \cong {\Biext}^1(A_1,A_2;Y_3(1)).$
More precisely we have the following isomorphisms
 \begin{equation}\label{geo1}
\begin{array}{c}
 {\Biext}^1(G_1,G_2;Y_3(1)) \cong {\Biext}^1(A_1,A_2;Y_3(1)) \\
 {\Biext}^1(G_1,G_2;G_3) \cong {\Biext}^1(G_1,G_2;Y_3(1))\\
\end{array}
\end{equation}
According to Definition~\ref{defmor} these isomorphisms mean that
 \emph{biextensions of 1-motives by 1-motives respect the weight filtration ${\W}_*$}, i.e. they satisfy the main property of morphisms of motives.

Inspired by~\cite{SGA7} Expos\'e VIII Corollary 2.2.11, if $M_1, M_2, M_3$ are three 1-motives defined over $S$, we require the anti-commutativity of the diagram
\[\begin{array}{ccc}
  {\Biext}^1(M_1,M_2;M_3)  & \stackrel{\cong}{\longrightarrow} & {\Biext}^1(M_2,M_1;M_3)  \\
{=} \downarrow  &  &  \downarrow {=} \\
 {\Hom}(M_1,M_2;M_3) & \stackrel{\cong}{\longrightarrow} & {\Hom}(M_2, M_1;M_3) \\
\end{array}\]
where the horizontal maps are induced by the morphism which permutes the factors.
The definitions of symmetric biextension of 1-motives~\ref{biextsym} and of skew-symmetric biextension of 1-motives~\ref{antisym} allow then the following definition:

\begin{defn}\label{defmorsym} Let $M$ and $M'$ be 1-motives defined over $S$.
\textbf{A symmetric morphism
$M \otimes M \rightarrow M'$} is an isomorphism class of skew-symmetric biextensions of $(M, M)$
by $M'$. \textbf{A skew-symmetric morphism
$M \otimes M \rightarrow M'$} is an isomorphism class of symmetric biextensions of $(M, M)$
by $M'$.
\end{defn}

Now we generalize Definition~\ref{defmor} to a finite tensor product of 1-motives:

\begin{lem}\label{otimesM}
For 1-motives defined over $S$, assume the existence of a weight filtration and of a tensor product, which are compatible one with another.
Let $l$ and $i$ be positive integers and let $M_j=[X_j \stackrel{u_j}{\longrightarrow} G_j]$ (for $j=1, \dots,l$) be a 1-motive defined over $S$. If $i \geq 1$ and $l +1 \geq i$, the motive $ \otimes^l_{j=1}M_j /{\W}_{-i}(\otimes^l_{j=1}M_j)$ is isogeneous to the motive \begin{equation}\label{formulaotimes}
  \sum
 \Big(\otimes_{k \in \{\nu_1, \dots,\nu_{l-i+1}\}}X_k \Big) \bigotimes
 \Big(\otimes_{j\in \{\iota_1, \dots, \iota_{i-1}\}} M_j /{\W}_{-i}(\otimes_{j\in \{\iota_1 , \dots, \iota_{i-1}\}}M_j) \Big)
\end{equation}
\pn where the sum is taken over all the $(l-i+1)$-uplets $\{\nu_1, \dots,\nu_{l-i+1}\}$ and all the $(i-1)$-uplets $\{\iota_1, \dots, \iota_{i-1}\}$ of $\{1, \cdots,l\}$ such that $ \{\nu_1, \dots,\nu_{l-i+1}\}\cap \{\iota_1, \dots, \iota_{i-1}\} = \emptyset$ and
$ \nu_1 < \dots < \nu_{l-i+1}$, $  \iota_1 < \dots < \iota_{i-1}.$
\end{lem}

\begin{proof}
1-motives $M_j$ have components of weight 0 (the lattice part $X_j$), of weight -1 (the abelian part $A_j$) and of weight -2 (the toric part $Y_j(1)$).
Consider the pure motive ${\Gr}^{\W}_{-i} (\otimes^l_{j=1}M_j)$: it is a finite sum of tensor products of $l$ factors of weight 0, -1 other -2.
If $i=l$ the tensor product
\[A_1 \otimes A_2 \otimes \dots \otimes A_l\]
\pn contains no factors of weight 0. For each $i$ strictly bigger than $l$, it is also easy to construct a tensor product of $l$ factors whose total weight is $-i$ and in which no factor has weight 0 (for example if
$i=l+2$ we take
\[Y_1(1) \otimes Y_2(1) \otimes A_3 \otimes \dots \otimes A_l).\]
\pn However if $i$ is strictly smaller than $l$,
in each of these tensor products of $l$ factors, there is at least one factor of weight 0, i.e. one of the $X_j$ for $j=1, \dots,l$.\\
Now fix a $i$ strictly smaller than $l$. The tensor products where there are less factors of weight 0 are exactly those where there are more factors of weight -1. Hence in the pure motive ${\Gr}_{-i} (\otimes^l_{j=1}M_j)$, the tensor products with less factors of weight 0 are of the type
\[ X_{\nu_1} \otimes \dots \otimes X_{\nu_{l-i}} \otimes A_{\iota_1} \otimes \dots \otimes A_{\iota_i}.\]
After these observations, the conclusion is clear. Remark that we have only an isogeny because in the 1-motive (\ref{formulaotimes}) the factor
\[X_{\nu_1} \otimes X_{\nu_2} \otimes \dots \otimes X_{\nu_p} \otimes
{\mathcal{Y}}_{\iota_1} \otimes {\mathcal{Y}}_{\iota_2} \otimes  \dots \otimes {\mathcal{Y}}_{\iota_{l-p}}\]
\pn appears with multiplicity ``$p+m$'' where $m$ is the number of
${\mathcal{Y}}_{\iota_q}$ (for $q=1, \dots,l-p$) which are of weight 0, instead of appearing only once like in the 1-motive \\
$ \otimes_{j}M_j /{\W}_{-i}(\otimes_{j}M_j)$.
In particular for each $i$ we have that
\begin{eqnarray}
\nonumber  {\Gr}^{\W}_0 \Big(
\sum (\otimes_{k }X_k ) \otimes (\otimes_{j} M_j /{\W}_{-i}) \Big) &=& l ~{\Gr}^{\W}_0 \Big( \otimes_{j}M_j /{\W}_{-i} \Big) \\
\nonumber  {\Gr}^{\W}_{-1} \Big(
\sum (\otimes_{k }X_k ) \otimes (\otimes_{j} M_j /{\W}_{-i}) \Big) &=& (l-1) ~{\Gr}^{\W}_{-1} \Big( \otimes_{j}M_j /{\W}_{-i} \Big)
\end{eqnarray}
\end{proof}

\begin{thm} \label{thmotimes}
For 1-motives defined over $S$, assume the existence of a weight filtration and of a tensor product, which are compatible one with another. Moreover, assume that the morphisms between 1-motives respect the weight filtration.
Let $M$ and $M_1, \dots, M_l$ be 1-motives over $S$.
Modulo isogenies a morphism from the tensor product of $M_1,\dots,M_l$
to $M$ is a sum of copies of isomorphism classes of biextensions of $(M_i, M_j)$ by $M$ for $i,j=1, \dots l$ and $i \not= j$.
More precisely we have that
\[{\Hom}(M_1,M_2, \ldots,M_l;M)\otimes {\QQ} =\sum {\Biext}^1(M_{\iota_1},M_{\iota_2};X_{\nu_1}^{\vee}\otimes \dots \otimes X_{\nu_{l-2}}^{\vee} \otimes M) \]
\pn where the sum is taken over all the $(l-2)$-uplets $\{\nu_1, \dots,\nu_{l-i+1}\}$ and all the 2-uplets $\{\iota_1,\iota_{2}\}$ of $\{1, \cdots,l\}$ such that $ \{\nu_1, \dots,\nu_{l-2}\}\cap \{\iota_1, \iota_{2}\} = \emptyset$ and
$ \nu_1 < \dots < \nu_{l-2}$, $  \iota_1 < \iota_{2}.$
\end{thm}

\begin{proof} Because morphisms of motives have to respect weights,
the only non trivial components of the morphism $\otimes^l_{j=1}M_j \rightarrow M$ are the components of the morphism
\[ \otimes^l_{j=1}M_j \Big/ {\W}_{-3}(\otimes^l_{j=1}M_j) \longrightarrow M.\]
Using the equality obtained in Lemma \ref{otimesM} with $i=-3$, we can write explicitly this last morphism in the following way
\[ \sum_{ \iota_1 < \iota_2 ~\mathrm{and}~ \nu_1< \dots < \nu_{l-2}
 \atop \iota_1 , \iota_2 \notin \{\nu_1, \dots,\nu_{l-2}\} }
 X_{\nu_1}\otimes \dots \otimes X_{\nu_{l-2}}\otimes ( M_{\iota_1} \otimes M_{\iota_2} /{\W}_{-3}( M_{\iota_1} \otimes M_{\iota_2})) \longrightarrow M.\]
To have the morphism
$$X_{\nu_1}\otimes \dots \otimes X_{\nu_{l-2}}\otimes ( M_{\iota_1} \otimes M_{\iota_2} /{\W}_{-3}( M_{\iota_1} \otimes M_{\iota_2})) \longrightarrow M$$
is equivalent to have the morphism
$$ M_{\iota_1} \otimes M_{\iota_2} /{\W}_{-3}( M_{\iota_1} \otimes M_{\iota_2}) \longrightarrow X_{\nu_1}^{\vee}\otimes \dots \otimes X_{\nu_{l-2}}^{\vee}\otimes M$$
where $X_{\nu_k}^{\vee}$ is the $S$-group scheme ${\uHom}( X_{\nu_k},{\ZZ})$ for $k=1, \ldots,l-2$. But as observed in 1.1 ``to tensor a motive by a motive of weight zero'' means to take a certain number of copies of this motive, and so from definition \ref{defmor} we get the expected conclusion.
\end{proof}

\subsection{Linear morphisms and pairings}

Let $A$ and $B$ be abelian $S$-schemes. According to~\cite{SGA7} Expos\'e VIII 3.2 we have the well-known canonical isomorphisms
\[{\Hom}(A,B) \cong {\Biext}^1(A,B^*;{\ZZ}(1)) \cong
{\Hom}(B^*,A^*)\]
where $A^*$ and $B^*$ are the
 Cartier duals of $A$ and $B$ respectively.
In the case where $B = A$,
through these canonical isomorphisms the Poincar\'e biextension of $A$,
denoted by ${\mathcal{P}}_{1,A}$,
corresponds to the identities morphisms $id_A:A \rightarrow A$
and $id_{A^*}:A^* \rightarrow A^*.$
More in general,
to a morphism $f:A \rightarrow B$ is associated
the pull-back $(f,id)^* {\mathcal{P}}_{1,B}$ via $(f,id)$ of the Poincar\'e biextension of $B$,
that we denote by ${\mathcal{P}}_{f,B}.$
To the transpose morphism $f^t:B^* \rightarrow A^*$ of $f$ is
 associated the pull-back via $(id, f^t)$ of ${\mathcal{P}}_{1,A}$,
that we denote by ${\mathcal{P}}_{f^t,A}.$ Clearly these two biextensions
are isomorphic:
\[{\mathcal{P}}_{f,B} \cong {\mathcal{P}}_{f^t,A}.\]

According to definition \ref{defmor}, a biextension of $(A,B^*)$
 by ${\ZZ}(1)$ is a morphism from $A \otimes B^*$ to ${\ZZ}(1)$:
\[ {\Hom}(A,B^*;{\ZZ}(1)) =
{\Biext}^1(A,B^*;{\ZZ}(1)). \]
In the case where $B = A$, the Poincar\'e biextension ${\mathcal{P}}_{1,A}$
 of $A$ is \textbf{the motivic Weil pairing
$A \otimes A^* \rightarrow {\ZZ}(1)$ of $A$}: we write it $e_{1,A}.$
The biextension ${\mathcal{P}}_{f,B}$ of $(A,B^*)$ by ${\ZZ}(1)$
is the pairing
\[e_{1,B} \circ (f \times id):A \otimes B^* \longrightarrow {\ZZ}(1).\]
We denote this pairing $ f \otimes B^*$. [The reason of this notation is that if we were in a Tannakian category, we could recover this pairing composing the morphism $f \times id$ with the evaluation morphism $ev_B:B \otimes B^{\vee} \longrightarrow 1$ of $B:$
\[ A \otimes B^* \stackrel{f \times id}{\longrightarrow} B \otimes B^* = B \otimes B^{\vee}  \otimes {\ZZ}(1)  \stackrel{ev_B \times id}{\longrightarrow} 1 \otimes  {\ZZ}(1) = {\ZZ}(1)].\]
In an analogous way, the biextension ${\mathcal{P}}_{f^t,A}$ is the pairing
$e_{1,A} \circ (id \times f^t)$ that we denote $A \otimes f^t$.
Since the biextensions ${\mathcal{P}}_{f,B}$ and $ {\mathcal{P}}_{f^t,A}$
are isomorphic we have that
\[ f \otimes B^* = A \otimes f^t.\]

\begin{lem}
Let $f: A \rightarrow B$ be a morphism of abelian $S$-schemes
 and let $f^t: B^* \rightarrow A^*$ be its transpose morphism.
The morphisms $f$ and $f^t$ are adjoint with respect to the motivic Weil Pairing. In particular, if the morphism $f$ has an inverse,
its inverse $f^{-1}:B \rightarrow A$ and its contragradient
$\widehat{f}=(f^{-1})^t: A^* \rightarrow B^*$ are adjoint for the motivic Weil Pairing.
\end{lem}

\begin{proof} The equality $f \otimes B^* = A \otimes f^t$ means that the following diagram is commutative
\[\begin{array}{ccc}
  A \otimes B^*  & \stackrel{id \times f^t}{\longrightarrow} &A \otimes A^*  \\
{\scriptstyle{f \times id}} \downarrow  &  &  \downarrow {\scriptstyle{e_{1,A}}} \\
  B \otimes B^* & \stackrel{e_{1,B}}{\longrightarrow} & {\ZZ}(1). \\
\end{array}\]
\end{proof}

If $S$ is the spectrum of a field $k$ of characteristic 0,
we can consider the Tannakian category $\langle A \rangle^\otimes$
generated by $A$ in an appropriate category of realizations. The motivic Galois group of $A$ is the motivic affine group scheme ${\Sp}(\Lambda),$ where $\Lambda$ is an element of $\langle A \rangle^\otimes$ endowed with the following universal property: for each object $X$ of $\langle A \rangle^\otimes$ there exists a morphism
$ X \rightarrow \Lambda \otimes X$ functorial in $X$ (see~\cite{D4} 8.4, 8.10, 8.11 (iii)). In the following proposition we discuss the main properties of the motivic Weil pairing $e_{1,A}$ of $A$ and in particular its link with the motivic Galois group of $A$.

\begin{prop} The motivic Weil pairing $e_{1,A}$ of $A$ is skew-symmetric and non-degenerate. Moreover, if $S$ is the spectrum of a field $k$ of characteristic 0, the motivic Weil pairing is invariant under the action of the motivic Galois group  of $A$.
\end{prop}

\begin{proof}
Since the Poincar\'e biextension of $A$ is a symmetric biextension,
by definition~\ref{defmorsym} the corresponding pairing is skew-symmetric.
The reason of the non-de\-ge\-ne\-ra\-cy of the pairing $e_{1,A}$ is that the Poincar\'e biextension ${\mathcal{P}}_{1,A}$ is trivial only if restricted to $A \times \{0\}$ and $ \{0\} \times A^*$. The pairing $e_{1,A}$ is an element of ${\Hom} (A, A^*;{\GG}_m)$, which can be viewed as an Artin motive since ${\Hom} (A, A^*;{\GG}_m)\cong A\otimes A^*\otimes {\GG}_m^{\vee}$ is of weight 0 (here ${\GG}_m^{\vee}$ is the Tannakian dual of ${\GG}_m$). Therefore the motivic Galois group of $A$ acts on ${\Hom} (A, A^*;{\GG}_m)$ via the Galois group $\gal$. Since the Poincar\'e biextension ${\mathcal{P}}_{1,A}$ is defined over $k$, also the corresponding pairing $e_{1,A}$ is defined over $k$, and therefore $e_{1,A}$ is invariant under the action of $\gal$, i.e. under the action of the motivic Galois group of $A$.
\end{proof}

Let $M_i=(X_i,Y^{\vee}_i,A_i,A^*_i, v_i ,v^*_i,\psi_i)= [X_i  \stackrel{u_i}{\longrightarrow}  G_i ]$ (for $i=1,2$) be a 1-motive over $S$.
According to~\cite{D1} (10.2.14) a morphism from $M_1$ to
$M_2$ is a 4-uplet of morphisms
$F=(f:A_1 \rightarrow A_2 ,~f^t:A_2^* \rightarrow A_1^*,~g:X_1 \rightarrow X_2,~h:Y^{\vee}_2 \rightarrow Y^{\vee}_1)$ where
\begin{itemize}
    \item $f$ is a morphism of abelian $S$-schemes with transpose morphism $f^t,$ and $g$ and $h$ are morphisms of character groups of $S$-tori;
    \item $f \circ v_1 =v_2 \circ g$ and dually $ f^t \circ v_2^*= v_1^* \circ h;$
    \item via the isomorphism ${\mathcal{P}}_{f^t,A_1}={\mathcal{P}}_{f,A_2},$ we have
$\psi_1(x_1,h(y^*_2))=\psi_2(g(x_1),y^*_2)$ for each $(x_1,y^*_2) \in X_1\times Y^{\vee}_2.$
\end{itemize}
The transpose morphism $F^t: M_2^* \rightarrow M_1^*$  of $F=(f,f^t,g,h)$ is
$(f^t:A_2^* \rightarrow A_1^*,~f:A_1 \rightarrow A_2 ,~h^{\vee}:Y_1 \rightarrow Y_2,~ g^{\vee}:X_2^{\vee} \rightarrow X_1^{\vee})$ where $ h^{\vee}$ and $g^{\vee}$ are the dual morphisms of $h$ and $g$, i.e. morphisms of cocharacter groups of $S$-tori.

As for abelian $S$-schemes, also for 1-motives
we have the following isomorphisms:
\[{\Hom}(M_1,M_2) \cong {\Hom}(M_1, [{\ZZ}\rightarrow 0];M_2)\cong {\Hom}(M_1, M_2^*;{\ZZ}(1)).\]
In fact using the definition of bilinear morphisms~\ref{defmor}, we prove that

\begin{prop}\label{linmorM-1}
Let $M_1$ and $M_2$ be two 1-motives defined over $S$.
\[ {\Hom}(M_1, M_2) \cong
{\Biext}^1(M_1,M_2^*;{\ZZ}(1)). \]
In other words, the biextensions of $(M_1, M_2^*)$
by ${\ZZ}(1)$ are the morphisms from $M_1$ to $M_2$.
\end{prop}

\begin{proof} A biextension of $(M_1, M_2^*)$
by ${\ZZ}(1)$ is $({\mathcal{P}}, \Gamma_1, \Gamma_2,0)$ where
${\mathcal{P}}$ is an object of ${\Biext}^1(A_1,A_2^*;{\ZZ}(1))$ and
 $\Gamma_1$ and $\Gamma_2$ are trivializations of the biextensions $(id_{A_1},v^*_2)^* {\mathcal{P}}$ and $(v_1,id_{A^*_2})^* {\mathcal{P}}$ respectively, which coincide over $X_1 \times Y_2^{\vee}.$
To have the biextension ${\mathcal{P}}$ is the same thing as to have a morphism
$f:A_1 \rightarrow A_2$ of abelian $S$-schemes with transpose morphism $f^t$. By~\cite{SGA7} Expos\'e VIII Proposition 3.7, to have the biextension $(v_1,id_{A^*_2})^* {\mathcal{P}}$ of $(X_1,A^*_2)$ by ${\ZZ}(1)$ (resp. $(id_{A_1},v^*_2)^* {\mathcal{P}}$ of $(A_1,Y_2^{\vee})$ by ${\ZZ}(1)$) is the same thing as to have a morphism $X_1 \rightarrow A_2$ (resp. $h:Y^{\vee}_2 \rightarrow A^*_1 $) equal to the composite $f \circ v_1 $ (resp. $ f^t \circ v_2^*$),
and this is the same thing as to have a morphism $g:X_1  \rightarrow X_2$ (resp. $h:Y^{\vee}_2 \rightarrow Y^{\vee}_1 $) such that $f \circ v_1 =v_2 \circ g$ (resp. $ f^t \circ v_2^*= v_1^* \circ h$). The condition that the two trivializations $\Gamma_1$ and $\Gamma_2$ coincide over $X_1 \times Y_2^{\vee}$ is equivalent to the condition $\psi_1(x_1,h(y^*_2))=\psi_2(g(x_1),y^*_2)$ for each $(x_1,y^*_2) \in X_1\times Y^{\vee}_2$.
\end{proof}

\begin{prop}\label{linmorM-2}
Let $M_1$ and $M_2$ be two 1-motives defined over $S$.
\[ {\Hom}(M_1, M_2) \cong
{\Biext}^1(M_1,[{\ZZ} \rightarrow 0];M_2). \]
In other words, the biextensions of $(M_1,[{\ZZ} \rightarrow 0])$
by $M_2$ are the morphisms from $M_1$ to $M_2$.
\end{prop}

\begin{proof}  A biextension $({\mathcal{P}},\Gamma_1, \Gamma_2,\lambda)$ of $(M_1,[{\ZZ} \rightarrow 0])$ by $M_2$ consists of a biextension of ${\mathcal{P}}$ of $(G_1, 0)$ by $G_2;$ a trivialization $\Gamma_1 $ (resp.  $\Gamma_2$) of the biextension $(u_1, id_{0})^*{\mathcal{P}}$ of $(X_1,0)$ by $G_2$  (resp. of the biextension
 $(id_{G_1},0)^*{\mathcal{P}}$ of $(G_1,{\ZZ})$ by $G_2$) such that $\Gamma_1 $ and $\Gamma_2 $ coincide over $X_1 \times {\ZZ}$; and a morphism $\lambda: X_1 \otimes {\ZZ} \rightarrow X_2$ such that the morphism $u_2 \circ \lambda: X_1 \otimes {\ZZ} \rightarrow G_2$ is compatible with the restriction $\Gamma$ of $\Gamma_1$ or $\Gamma_2$ over $X_1 \times {\ZZ}$, i.e. the following diagram is commutative
\begin{equation}\label{linmorM-2:1}
\begin{array}{ccccc}
  G_2 &= &  G_2 ~~~~~& & \\
{\scriptstyle \Gamma_2}  \uparrow ~~~ &  & {\scriptstyle \Gamma}  \uparrow  ~~~~~~ & \nwarrow  {\scriptstyle u_2} & \\
G_1 \times {\ZZ} & \stackrel{(u_1,id_{\ZZ})}{\longleftarrow} & X_1 \otimes {\ZZ} \stackrel{\lambda}{\longrightarrow} & X_2. \\
\end{array}
\end{equation}
The trivialization $\Gamma_2$ defines a morphism $\gamma$ from $G_1$ to $G_2$, the morphism $\lambda$ defines a morphism, again called $\lambda$, from $X_1$ to $X_2$ and the commutativity of the above diagram implies the commutativity of the diagram
\[
\begin{array}{ccc}
 X_1 & \stackrel{u_1}{\longrightarrow} & G_1  \\
 {\scriptstyle \lambda} \downarrow   &  & \downarrow {\scriptstyle \gamma} \\
 X_2 & \stackrel{u_2}{\longrightarrow} & G_2. \\
\end{array}
\]
\end{proof}

Recall also that to each 1-motive $M=(X_1,Y^{\vee}_1,A_1,A^*_1, v_1 ,v^*_1,\psi_1)=[X \stackrel{u}{\longrightarrow }G]$ is associated its Poincar\'e biextension, that we denote by ${\mathcal{P}}_{1,M}$, which expresses the Cartier duality between $M$ and $M^*$. It is the biextension
$({\mathcal{P}}_{1,A}, \psi_1, \psi_2,0)$ of $(M,M^*)$ by ${\ZZ}(1)$ where
 $\psi_1$ is the trivialization of the biextension $(id_A,v^*)^* {\mathcal{P}}_{1,A}$ which
defines the morphism $u: X \rightarrow G$, and $\psi_2$ is the trivialization of the biextension $(v,id_{A^*})^* {\mathcal{P}}_{1,A}$ which defines the morphism $u^*: Y^{\vee} \rightarrow G^*.$

Via the isomorphism of Proposition~\ref{linmorM-1} the Poincar\'e biextension of $M_1$,  ${\mathcal{P}}_{1,M_1}$,
corresponds to the identities morphisms $id_{M_1}:M_1 \rightarrow M_1$
and $id_{M_1^*}:M_1^* \rightarrow M_1^*.$
More in general,
to a morphism $F=(f,f^t,g,h):M_1 \rightarrow M_2$ is associated
the pull-back $(F\times id)^* {\mathcal{P}}_{1,M_2}$ by $F\times id$ of the Poincar\'e biextension of $M_2$,
that we denote by ${\mathcal{P}}_{F,M_2}.$
Explicitly, if $({\mathcal{P}}_{1,A_2}, \psi_1^2, \psi_2^2,0)$ is the Poincar\'e biextension of $M_2$, the biextension ${\mathcal{P}}_{F,M_2}$ is
\[((f\times id)^*{\mathcal{P}}_{1,A_2},(f\times id)^* \psi_1^2, (g\times id)^*\psi_2^2,0).\]
 To the transpose morphism $F^t=(f^t,f,h^*,g^*):M_2^* \rightarrow M_1^*$ of $F$ is
 associated the pull-back via $id\times F^t$ of ${\mathcal{P}}_{1,M_1}$,
that we denote by ${\mathcal{P}}_{F^t,M_1}.$
Explicitly, if $({\mathcal{P}}_{1,A_1}, \psi_1^1, \psi_2^1,0)$ is the Poincar\'e biextension of $M_1$, the biextension ${\mathcal{P}}_{F^t,M_1}$ is
\[((id\times f^t)^*{\mathcal{P}}_{1,A_1}, (id\times h)^* \psi_1^1, (id\times f^t)^*\psi_2^1,0).\]
As for abelian schemes we have
\[{\mathcal{P}}_{F,M_2} \cong {\mathcal{P}}_{F^t,M_1}.\]

According to definition \ref{defmor}, each biextensions of $(M_1,M_2^*)$
 by ${\ZZ}(1)$ is a morphism from $M_1 \otimes M_2^*$ to ${\ZZ}(1)$:
 \[ {\Hom}(M_1, M_2^*;{\ZZ}(1)) =
{\Biext}^1(M_1,M_2^*;{\ZZ}(1)). \]
In the case where $M_1 = M_2$, the Poincar\'e biextension ${\mathcal{P}}_{1,M_1}$
 of $M_1$ is \textbf{the motivic Weil pairing
$M_1 \otimes M_1^* \rightarrow {\ZZ}(1)$ of $M_1$}: we write it $e_{1,M_1}.$
The biextension ${\mathcal{P}}_{F,M_2}$ of $(M_1,M_2^*)$ by ${\ZZ}(1)$
is the pairing
\[e_{1,M_2} \circ (F \times id):M_1 \otimes M_2^* \longrightarrow {\ZZ}(1).\]
We denote this pairing $ F \otimes M_2^*$. \\
In an analogous way, the biextension ${\mathcal{P}}_{F^t,M_1}$ is the pairing
$e_{1,M_1} \circ (id \times F^t)= M_1 \otimes F^t$.
Since the biextensions ${\mathcal{P}}_{F,M_2}$ and $ {\mathcal{P}}_{F^t,M_1}$
are isomorphic we have that
\[ F \otimes M_2^* = M_1 \otimes F^t.\]
As for abelian schemes this last equality implies

\begin{lem}
Let $F: M_1 \rightarrow M_2$ be a morphism of abelian $S$-schemes
 and let $F^t: M_2^* \rightarrow M_1^*$ be its transpose morphism.
The morphisms $F$ and $F^t$ are adjoint with respect to the motivic Weil Pairing.
 In particular, if the morphism $F$ has an inverse,
its inverse $F^{-1}:M_2 \rightarrow M_1$ and its contragradient
$\widehat{F}=(F^{-1})^t: M_1^* \rightarrow M_2^*$ are adjoint for the motivic Weil Pairing.
\end{lem}

\section{Realizations of biextensions}

\subsection{Construction of the Hodge realization of biextensions}

Let $S$ be the spectrum of the field $\CC$ of complex numbers.
Recall that a \textbf{mixed Hodge structure $({\h}_{\ZZ}, {\W}_*, {\F}^*)$} consists of a finitely generated $\ZZ$-module ${\h}_\ZZ$, an increasing filtration $W_*$ (the weight filtration) on ${\h}_{\ZZ} \otimes \QQ$, a decreasing filtration ${\F}^*$ (the Hodge filtration) on ${\h}_{\ZZ} \otimes \CC$, and  some axioms relating these two filtrations (see~\cite{D71} D\'efinition 1.1).
Let $M_i=(X_i,A_i,Y_i(1),G_i,u_i)$ (for $i=1,2,3$) be a 1-motive over $\CC$. The Hodge realization $T_{\h}(M_i)=(T_{\Bbb Z}(M_i), {\W}_*,{\F}^*)$ of the 1-motive $M_i$ is the mixed Hodge structure consisting of the fibred product
$T_{\Bbb Z}(M_i)={\rm Lie}(G_i)\times_{G_i} X_i $ (viewing ${\rm Lie}(G_i)$ over $G_i$ via the exponential map and $X_i$ over $G_i$ via $u_i$) and of the weight and Hodge filtrations defined in the following way:
\begin{eqnarray}
\nonumber  {\W}_{0}(T_{\Bbb Z}(M_i)) &=& T_{\Bbb Z}(M_i), \\
\nonumber  {\W}_{-1}(T_{\Bbb Z}(M_i)) &=& {\h}(G_i,\ZZ), \\
\nonumber  {\W}_{-2}(T_{\Bbb Z}(M_i)) &=& {\h}( Y_i(1),\ZZ), \\
\nonumber  {\F}^0(T_{\Bbb Z}(M_i) \otimes {\CC}) &=& \ker( T_{\Bbb Z}(M_i) \otimes {\CC} \longrightarrow {\rm Lie}(G_i)).
\end{eqnarray}
(see~\cite{D1} \S 10.1.3). Denote by $M_i^{\an}$ the complex of analytic groups $[X_i \rightarrow G_i^{\an}].$
Each biextension of $(M_1,M_2)$ by $M_3$ defines a unique biextension of $(M_1^{\an},M_2^{\an})$ by $M_3^{\an}$:

\begin{prop}\label{We}
Let $M_i=(X_i,A_i,Y_i(1),G_i,u_i)$ (for $i=1,2,3$) be a 1-motive over $\CC$.
 The application
\[{\Biext}^1(M_1,M_2;M_3) \longrightarrow {\Biext}^1(M_1^{\an},M_2^{\an};M_3^{\an})\]
 is injective. Its image consists of the biextensions $({\mathcal{B}}, \Psi_1, \Psi_2,\lambda)$ whose restriction ${\mathcal{B}}$ to $G_1 \times G_2$ comes, via pull-backs and push-downs, from a biextension $B$ of $(A_1,A_2)$ by $Y_3(1)$:
\[{\mathcal{B}} = \iota_{3\,*}(\pi_1,\pi_2)^*B\]
 where $\pi_i:G_i \rightarrow A_i$
is the projection of $G_i$ over $A_i$ (for $i=1,2$) and $\iota_3:Y_3(1) \rightarrow
G_3$ is the inclusion of $Y_3(1)$ over $G_3.$
\end{prop}

The morphism ${\T}_{\ZZ}(M_1)\otimes {\T}_{\ZZ}(M_2)\rightarrow  {\T}_{\ZZ}(M_3)$ corresponding to a biextension $({\mathcal{B}}, \Psi_1, \Psi_2,\lambda)$ of $(M_1,M_2)$ by $M_3$ comes from the trivializations defining the biextension of $(M_1^{\an},M_2^{\an})$ by $M_3^{\an}$ induced by $({\mathcal{B}}, \Psi_1, \Psi_2,\lambda)$. Therefore in order to find the Hodge realization of a biextension $({\mathcal{B}}, \Psi_1, \Psi_2,\lambda)$ of $(M_1,M_2)$ by $M_3$ we first have to compute ${\Biext}^1(M_1^{\an},M_2^{\an};M_3^{\an}):$

\begin{lem}\label{Ho}
Let $M_i$ (for $i=1,2,3$) be a 1-motive over $\CC$.
The group\\ ${\Biext}^1(M_1^{\an},M_2^{\an};M_3^{\an})$ is isomorphic to the group of the applications
\[ \Phi:{\T}_{\ZZ}(M_1)\otimes {\T}_{\ZZ}(M_2) \longrightarrow {\T}_{\ZZ}(M_3)\]
such that $\Phi_{\CC}: {\T}_{\CC}(M_1)\otimes {\T}_{\CC}(M_2) \rightarrow {\T}_{\CC}(M_3)$ respects the Hodge filtration.
\end{lem}

\begin{proof} The proof will be done in several steps.\\
 \textbf{Step 1:} Let $V_{\CC}, W_{\CC},Z_{\CC}$ be three vector spaces. Since extensions of vector spaces are trivial, by \cite{SGA7} Expos\'e VIII 1.5 we have that for $i=0,1$
\begin{eqnarray}
\nonumber   {\Biext}^i (V_{\CC},W_{\CC};Z_{\CC}) &=& {\Ext}^i ( V_{\CC},{\uHom}(W_{\CC},Z_{\CC}))\\
\nonumber    &=& {\Ext}^i(V_{\CC},W_{\CC}^* \otimes Z_{\CC}).
\end{eqnarray}
Therefore in the analytic category we obtain that
\begin{eqnarray}
 \label{step1-1} {\Biext}^0 (V_{\CC},W_{\CC};Z_{\CC}) & \cong & {\Hom}(V_{\CC} \otimes W_{\CC}, Z_{\CC}) \\
 \label{step1-2}  {\Biext}^1 (V_{\CC},W_{\CC};Z_{\CC}) &=& 0.
\end{eqnarray}
 Let $V_{\ZZ},W_{\ZZ},Z_{\ZZ}$ be three free finitely generated $\ZZ$-modules contained respectively in $V_{\CC},W_{\CC},Z_{\CC}$. Since the morphism of complexes
$ [V_{\ZZ} \rightarrow V_{\CC}] \rightarrow  [0 \rightarrow V_{\CC}/ V_{\ZZ}] $ is a quasi-isomorphism one can check that we have the equivalence of categories
\begin{equation}\label{qiso}
{\bBiext}(V_{\CC}/ V_{\ZZ},W_{\CC}/ W_{\ZZ};Z_{\CC}/ Z_{\ZZ}) \cong
{\bBiext}([V_{\ZZ} \rightarrow V_{\CC}],[W_{\ZZ} \rightarrow W_{\CC}];[Z_{\ZZ} \rightarrow Z_{\CC}])
\end{equation}
In order to get such an equivalence, one can also use the homological interpretation of biextensions stated in Remark~\ref{homolo} .
Now we will prove that
\begin{equation}\label{base}
{\Biext}^1(V_{\CC}/ V_{\ZZ},W_{\CC}/ W_{\ZZ};Z_{\CC}/ Z_{\ZZ}) \cong
{\Hom}(V_{\ZZ} \otimes W_{\ZZ}, Z_{\ZZ}).
\end{equation}
Let $\phi_i: V_{\CC} \otimes W_{\CC} \rightarrow Z_{\CC}$ ($i=1,2$)
be a bilinear application such that the restriction of $\phi_1 -\phi_2$ to
$V_{\ZZ} \otimes W_{\ZZ}$ factors through $Z_{\ZZ} \rightarrow Z_{\CC},$ i.e. it takes values in $Z_{\ZZ}$. Denote by ${\mathcal{B}}(\phi_1,\phi_2)$ the following biextension of $([V_{\ZZ} \rightarrow V_{\CC}],[W_{\ZZ} \rightarrow W_{\CC}])$ by $[Z_{\ZZ} \rightarrow Z_{\CC}]$:
the trivial biextension $V_{\CC} \times W_{\CC} \times Z_{\CC}$ of
$(V_{\CC},W_{\CC})$ by $Z_{\CC}$, its trivializations
$\phi_1: V_{\ZZ} \otimes W_{\CC} \rightarrow Z_{\CC}$ and $\phi_2: V_{\CC} \otimes W_{\ZZ} \rightarrow Z_{\CC}$ and the morphism $\Phi: V_{\ZZ} \times W_{\ZZ} \rightarrow Z_{\ZZ}$ compatible with the trivializations $\phi_i$ (for $i=1,2$), i.e.
$\Phi=\phi_1-\phi_2$. According to (\ref{step1-2}) each biextension is like that and by (\ref{step1-1})
 two biextensions ${\mathcal{B}}(\phi_1,\phi_2)$ and
${\mathcal{B}}(\phi'_1,\phi'_2)$ are isomorphic if and only if $\phi_1-\phi_2= \phi'_1-\phi'_2$.

\pn \textbf{Step 2:} Let $F_1,F_2,F_3$ be three subspaces of $V_{\CC},W_{\CC},Z_{\CC}$ respectively. Consider the complexes
\begin{center}
$ K_1= [V_{\ZZ}\oplus F_1 \longrightarrow V_{\CC}] \qquad \qquad K'_1= [V_{\ZZ} \longrightarrow V_{\CC}] $\\
$  K_2= [W_{\ZZ}\oplus F_2 \longrightarrow W_{\CC}] \qquad \qquad  K'_2= [W_{\ZZ} \longrightarrow W_{\CC}] $\\
$  K_3= [Z_{\ZZ}\oplus F_3 \longrightarrow Z_{\CC}] \qquad \qquad K'_3= [Z_{\ZZ} \longrightarrow Z_{\CC}] $
\end{center}
In this step we prove that to have a biextension of $(K_1,K_2)$ by $K_3$ is the same thing as to have:
\begin{description}
\item[$(a)$] a biextension $B$ of $(K'_1,K'_2)$ by $K'_3$;
\item[$(b)$] a unique determined trivialization $\phi_1$ (resp. $\phi_2$) of the biextension
of $([0 \rightarrow F_1], K'_2)$ (resp. $(K'_1, [0 \rightarrow F_2])$) by $K'_3,$ pull-back of $B$;
\item[$(c)$] a unique determined trivialization of the biextension of $(F_1,F_2)$ by $F_3$ whose push-down via the inclusion
$F_3 \rightarrow B_3$ coincides with the restriction of $\phi_1-\phi_2$ to $F_1 \times F_2.$
\end{description}
We start observing that for ($i=0,1$):
\begin{enumerate}
 \item ${\Biext}^i (F_1,K'_2;K'_3)={\Biext}^i (K'_1,F_2;K'_3)=0$:
for $i=0$ this is a consequence of the fact that a bilinear application $f:F_1 \times W_{\CC} \rightarrow Z_{\CC}$ such that $f(F_1,W_{\ZZ}) \subseteq Z_{\ZZ}$ is trivial. For the assertion with $i=1$ we use (\ref{step1-2}) and the fact that each biadditif morphism $F_1 \times W_{\ZZ} \rightarrow Z_{\CC}$ comes from a biadditif morphism $F_1 \times W_{\CC} \rightarrow Z_{\CC}$, i.e. the trivialization over $F_1 \times W_{\ZZ}$ lifts to $F_1 \times W_{\CC}$.    \item ${\Biext}^i (F_1,F_2;K'_3)=0$ for:
since the biextensions of $(F_1,F_2)$ by $K'_3$ are the restriction to $F_1\times F_2$
of the biextensions of $(F_1,K'_2)$ and of $(K'_1,F_2)$ by $K'_3$, we can conclude using (1).
\item ${\Biext}^i (K'_1,K'_2;F_3)=0$: for $i=0$ this is a consequence of the fact that a bilinear application $f:V_{\CC} \times W_{\CC} \rightarrow F_3$ such that $f(V_{\ZZ},W_{\ZZ}) = 0$ is trivial. The proof of the assertion with $i=1$ is the same as in (1).
\item ${\Biext}^i (F_1,K'_2;F_3)={\Biext}^i (K'_1,F_2;F_3)=0$:
the biextensions of $(F_1,K'_2)$ (resp. of $(K'_1,F_2)$) by $F_3$ are the restriction to $F_1\times K'_2$ (resp. to $K'_1\times F_2$)
of the biextensions of $(K'_1,K'_2)$ by $F_3$, and so we can conclude using (3).
\item ${\Biext}^i (F_1,K'_2;K_3)={\Biext}^i (K'_1,F_2;K_3)=0$:
these results follow from (1), (4) and from the long exact sequence
\[  \begin{array}{c}
 0 \rightarrow {\Biext}^0 (F_1,K'_2;F_3)\rightarrow {\Biext}^0 (F_1,K'_2;K_3)\rightarrow {\Biext}^0 (F_1,K'_2;K'_3) \rightarrow \\
\rightarrow {\Biext}^1 (F_1,K'_2;F_3)\rightarrow {\Biext}^1 (F_1,K'_2;K_3)\rightarrow
{\Biext}^1(F_1,K'_2;K'_3) \rightarrow ...\\
\end{array}\]
\end{enumerate}
Using the exact sequences $0 \rightarrow F_i \rightarrow K_i \rightarrow K'_i \rightarrow 0$ (for $i=1,2,3$), we have the long exact sequences
\begin{equation}\label{step2-1}
  \begin{array}{c}
 0 \rightarrow {\Biext}^0 (K'_1,K_2;K_3)\rightarrow {\Biext}^0 (K_1,K_2;K_3)\rightarrow {\Biext}^0 (F_1,K_2;K_3) \rightarrow \\
\rightarrow {\Biext}^1 (K'_1,K_2;K_3)\rightarrow {\Biext}^1 (K_1,K_2;K_3)\rightarrow
{\Biext}^1(F_1,K_2;K_3) \rightarrow ...\\
\end{array}
\end{equation}
\begin{equation}\label{step2-2}
  \begin{array}{c}
 0 \rightarrow {\Biext}^0 (F_1,K'_2;K_3)\rightarrow {\Biext}^0 (F_1,K_2;K_3)\rightarrow {\Biext}^0 (F_1,F_2;K_3) \rightarrow \\
\rightarrow {\Biext}^1 (F_1,K'_2;K_3)\rightarrow {\Biext}^1 (F_1,K_2;K_3)\rightarrow
{\Biext}^1(F_1,F_2;K_3) \rightarrow ...\\
\end{array}
\end{equation}
\begin{equation}\label{step2-3}
  \begin{array}{c}
 0 \rightarrow {\Biext}^0 (F_1,F_2;F_3)\rightarrow {\Biext}^0 (F_1,F_2;K_3)\rightarrow {\Biext}^0 (F_1,F_2;K'_3) \rightarrow \\
\rightarrow {\Biext}^1 (F_1,F_2;F_3)\rightarrow {\Biext}^1 (F_1,F_2;K_3)\rightarrow
{\Biext}^1(F_1,F_2;K'_3) \rightarrow ...\\
\end{array}
\end{equation}
\begin{equation}\label{step2-4}
  \begin{array}{c}
 0 \rightarrow {\Biext}^0 (K'_1,K'_2;K_3)\rightarrow {\Biext}^0 (K'_1,K_2;K_3)\rightarrow {\Biext}^0 (K'_1,F_2;K_3) \rightarrow \\
\rightarrow {\Biext}^1 (K'_1,K'_2;K_3)\rightarrow {\Biext}^1 (K'_1,K_2;K_3)\rightarrow
{\Biext}^1(K'_1,F_2;K_3) \rightarrow ...\\
\end{array}
\end{equation}
\begin{equation}\label{step2-5}
  \begin{array}{c}
 0 \rightarrow {\Biext}^0 (K'_1,K'_2;F_3)\rightarrow {\Biext}^0 (K'_1,K'_2;K_3)\rightarrow {\Biext}^0 (K'_1,K'_2;K'_3) \rightarrow \\
\rightarrow {\Biext}^1 (K'_1,K'_2;F_3)\rightarrow {\Biext}^1 (K'_1,K'_2;K_3)\rightarrow
{\Biext}^1(K'_1,K'_2;K'_3) \rightarrow ...\\
\end{array}
\end{equation}
From (3), (5), \ref{step2-4} and \ref{step2-5} (resp. (2), (5), \ref{step2-2} and \ref{step2-3}) we obtain the inclusions of categories
\begin{eqnarray}
\nonumber  {\bBiext} (K'_1,K_2;K_3) &\subseteq & {\bBiext} (K'_1,K'_2;K'_3)\\
\nonumber ({\mathrm{resp.}} \quad {\bBiext} (F_1,K_2;K_3)   &\subseteq& {\bBiext} (F_1,F_2;F_3)).
\end{eqnarray}
 Using (\ref{step2-1}) we can conclude.\\
According to step 1, we can reformulate what we have proved in the following way: the group
 ${\Biext}^1 (K_1,K_2;K_3)$ is isomorphic to the group of applications
\[ \Phi: V_{\ZZ}\otimes W_{\ZZ} \longrightarrow Z_{\ZZ}\]
such that $\Phi_{\CC}: V_{\CC} \otimes W_{\CC} \rightarrow Z_{\CC}$
satisfies $\Phi_{\CC}(F_1,F_2) \subseteq F_3$. Explicitly, the biextension
 of $(K_1,K_2)$ by $K_3$ associated via this
isomorphism to the application $\Phi: V_{\ZZ}\otimes W_{\ZZ} \rightarrow Z_{\ZZ},$
is the following one:
by step 1, a biextension of $(K'_1,K'_2)$ by $K'_3$ consists of
the trivial biextension $V_{\CC} \times W_{\CC} \times Z_{\CC}$ of
$(V_{\CC},W_{\CC})$ by $Z_{\CC}$, two of its trivializations
$\phi_1: V_{\ZZ} \otimes W_{\CC} \rightarrow Z_{\CC}$ and $\phi_2: V_{\CC} \otimes W_{\ZZ} \rightarrow Z_{\CC}$, where
$\phi_i: V_{\CC} \otimes W_{\CC} \rightarrow Z_{\CC}$ ($i=1,2$) is a
bilinear application, and a morphism $\Phi:V_{\ZZ} \otimes W_{\ZZ} \rightarrow Z_{\ZZ}$ compatible with the trivializations $\phi_i,$ i.e.
 $\Phi=\phi_1 -\phi_2$. According to step 2, this biextension of $(K'_1,K'_2)$ by $K'_3$ comes from a biextension of $(K_1,K_2)$ by $K_3$ if $\phi_2: F_1 \otimes W_{\CC} \rightarrow Z_{\CC}$ and $\phi_1: V_{\CC} \otimes F_2 \rightarrow Z_{\CC}$ are such that
$\Phi_{\CC}=\phi_1 -\phi_2(F_1,F_2) \subseteq F_3.$\\
In other words the biextension $({\mathcal{B}}, \Psi_1, \Psi_2,\lambda)$ of $(K_1,K_2)$ by $K_3$ associated to the application $\Phi: V_{\ZZ}\otimes W_{\ZZ} \rightarrow Z_{\ZZ}$ is defined in the following way:
the trivial biextension $V_{\CC} \times W_{\CC} \times Z_{\CC}$ of
$(V_{\CC},W_{\CC})$ by $Z_{\CC}$, its trivializations
\begin{equation}\label{hl}
  \begin{array}{c}
 \Psi_1: \big(V_{\ZZ} \oplus F_1\big) \times W_{\CC} \longrightarrow Z_{\CC}, ~~~~~~
(v_{\ZZ} \oplus f_1, w_{\CC})\mapsto \phi_1 (v_{\ZZ}, w_{\CC})+  \phi_2( f_1,w_{\CC}) \\
\Psi_2: V_{\CC} \times \big( W_{\ZZ} \oplus F_2 \big) \longrightarrow Z_{\CC},~~~~~~
(v_{\CC}, w_{\ZZ}\oplus f_2)\mapsto \phi_2 (v_{\CC},w_{\ZZ}) + \phi_1(v_{ \CC},f_2 )\\
\end{array}
\end{equation}
and the morphism $\lambda=\phi_1 -\phi_2: V_{\ZZ} \otimes W_{\ZZ} \rightarrow Z_{\ZZ}.$

\pn \textbf{Step 3:} In order to conclude we apply what we have proved in step 2 to the complexes $[ {\T}_{\ZZ}(M_i) \oplus {\F}^0{\T}_{\CC}(M_i) \rightarrow
{\T}_{\CC}(M_i)]$ (for $i=1,2,3$): in fact
\begin{itemize}
    \item for each 1-motive $M_i$ we have the quasi-isomorphisms
\[\begin{array}{ccc}
 {\T}_{\ZZ}(M_i) \oplus {\F}^0{\T}_{\CC}(M_i) & \longrightarrow & X_i  \\
 \downarrow  &  &  \downarrow \\
 {\T}_{\CC}(M_i)  & \longrightarrow & G_i; \\
\end{array}\]
    \item the only non trivial condition to check in order to prove that $\Phi_{\CC}: {\T}_{\CC}(M_1) \otimes {\T}_{\CC}(M_2) \rightarrow {\T}_{\CC}(M_3)$
respect the Hodge filtration $\F^*$ is
\[\Phi_{\CC} \big({\F}^0{\T}_{\CC}(M_1) \otimes {\F}^0{\T}_{\CC}(M_2)\big) \subseteq {\F}^0{\T}_{\CC}(M_3).\]
\end{itemize}
\end{proof}

 \emph{ Proof of Proposition~\ref{We}}. Recall that by (G.A.G.A)
\begin{equation} \label{gaga}
 {\bBiext}^1(A_1,A_2;Y_3(1)) \cong {\bBiext}^1(A_1^{\an},A_2^{\an};Y_3^{\an}(1)).
\end{equation}
We first prove the injectivity.
Let $({\mathcal{B}}, \Psi_1, \Psi_2,\lambda)$ be a biextension of $(M_1,M_2)$ by $M_3$ and let $B$ be the biextension of $(A_1,A_2)$ by $Y_3(1)$ corresponding to ${\mathcal{B}}$ via the equivalence of categories described in~\cite{B3} Theorem 2.5.2. Suppose that
$({\mathcal{B}}, \Psi_1, \Psi_2,\lambda)^{\an}$ is the trivial biextension
of $(M_1^{\an},M_2^{\an})$ by $M_3^{\an}$. According to
(\ref{gaga}) the biextension $B$ is trivial, and so because of~\cite{B3} Theorem 2.5.2 also the biextension ${\mathcal{B}}$ is trivial.
Hence the biextension $({\mathcal{B}}, \Psi_1, \Psi_2,\lambda)$ is defined through the biadditive applications $\Psi_1: X_1 \times G_2 \rightarrow G_3,~
 \Psi_2: G_1 \times X_2 \rightarrow G_3 $ and $\lambda: X_1 \times X_2 \rightarrow X_3.$ By hypothesis these applications are zero in the analytic category, and therefore they are zero. This prove the injectivity. \\
 Now let $({\mathcal{B}}, \Psi_1, \Psi_2,\lambda)$ be a biextension of $(M_1^{\an},M_2^{\an})$ by $M_3^{\an}$ satisfying the condition of this lemma.
We have to prove that it is algebraic.
Clearly the application $\lambda: X_1 \times X_2 \rightarrow X_3$ is algebraic.
 By (\ref{gaga}) and the equivalence of categories described in~\cite{B3} Theorem 2.5.2, the biextension ${\mathcal{B}}$ of $(G_1,G_2)$ by $G_3$ is algebraic.
In order to conclude we have to prove that also the trivializations $\Psi_1: X_1 \times G_2 \rightarrow G_3$ and $\Psi_2: G_1 \times X_2 \rightarrow G_3 $
of ${\mathcal{B}}$ are algebraic. But this is again a consequence of (G.A.G.A).$\Box$\\

Denote by $\mathcal{MHS}$ the category of mixed Hodge structures. Recall that a morphism $(H_{\Bbb Z}, {\W}_*,{\F}^*) \rightarrow (H_{\Bbb Z}', {\W}_*,{\F}^*)$ of mixed Hodge structures consists of a morphism
$f_{\ZZ}: H_{\Bbb Z} \rightarrow H_{\Bbb Z}'$ such that
$f_{\QQ}: H_{\Bbb Z} \otimes \, \QQ \rightarrow H_{\Bbb Z}' \otimes \, \QQ$
and $f_{\CC}: H_{\Bbb Z} \otimes \, \CC \rightarrow H_{\Bbb Z}' \otimes \, \CC$
are compatible with the weight filtration ${\W}_*$ and the Hodge filtration ${\F}^*$ respectively.

\begin{thm}\label{hodge}
Let $M_i$ (for $i=1,2,3$) be a 1-motive over $\CC$ and let $T_{\h}(M_i)=(T_{\Bbb Z}(M_i), {\W}_*,{\F}^*)$ be its Hodge realization.
We have that
\[{\Biext}^1(M_1,M_2;M_3) \cong {\Hom}_{\mathcal{MHS}}\big(
{\T}_{\h}(M_1)\otimes {\T}_{\h}(M_2),  {\T}_{\h}(M_3)\big).\]
\end{thm}

\begin{proof} By Lemma~\ref{Ho} we can identify the elements of
 ${\bBiext}^1(M_1^{\an},M_2^{\an};M_3^{\an})$ with applications
 $\Phi: {\T}_{\ZZ}(M_1)\otimes {\T}_{\ZZ}(M_2) \rightarrow {\T}_{\ZZ}(M_3)$ such that $\Phi_{\CC}: {\T}_{\CC}(M_1)\otimes {\T}_{\CC}(M_2) \rightarrow {\T}_{\CC}(M_3)$ is compatible with the Hodge filtration $\F^*$. Then Proposition~\ref{We} furnishes a bijection between ${\Biext}^1(M_1,M_2;M_3)$ and the set $\mathcal{H}$ of applications
\[ \Phi: {\T}_{\ZZ}(M_1)\otimes {\T}_{\ZZ}(M_2) \longrightarrow {\T}_{\ZZ}(M_3)\]
having the following properties
\begin{description}
    \item[$(a)$] $\Phi_{\CC}: {\T}_{\CC}(M_1)\otimes {\T}_{\CC}(M_2) \rightarrow {\T}_{\CC}(M_3)$ is compatible with the Hodge filtration $\F^*$;
    \item[$(b)$] the restriction of $\Phi$ to ${\W}_{-1}({\T}_{\ZZ}(M_1))\otimes {\W}_{-1}({\T}_{\ZZ}(M_2)) \rightarrow {\W}_{-1}({\T}_{\ZZ}(M_3))$ comes from a morphism ${\Gr}_{-1}({\T}_{\ZZ}(M_1))\otimes {\Gr}_{-1}({\T}_{\ZZ}(M_2)) \rightarrow {\Gr}_{-2}({\T}_{\ZZ}(M_3))$
 i.e. $\Phi$ is compatible with the weight filtration $\W_*.$
\end{description}
But by definition of morphisms in the category $\mathcal{MHS}$, the set $\mathcal{H}$ is nothing else as the group ${\Hom}_{\mathcal{MHS}}\big(
{\T}_{\h}(M_1)\otimes {\T}_{\h}(M_2),  {\T}_{\h}(M_3)\big).$
\end{proof}

\subsection{Construction of the $\ell$-adic realization of biextensions}

Let $S$ be the spectrum of a field $k$ of characteristic 0 embeddable in $\CC$. Let $M_i=(X_i,A_i,Y_i(1),$\\
$G_i,u_i)$ (for $i=1,2,3$) be a 1-motive over $k$.
We write it as a complex $[X_i {\buildrel u_i \over \longrightarrow} G_i]$ concentrated in degree 0 and 1. For each integer $n \geq 1$, let $[{\ZZ} {\buildrel n \over \longrightarrow} {\ZZ}]$ be the complex concentrated in degree -1 and 0. Consider the ${\ZZ}/n{\ZZ}\,$-module
\begin{eqnarray}
\nonumber  {\T}_{{\ZZ}/n{\ZZ}}(M_i) &=& {\h}^0(M_i {\otimes}^{\LL} {\ZZ}/n{\ZZ}) \\
\nonumber   &=& \bigg\{(x,g) \in X_i \times G_i ~~\vert~~ u_i(x)=n \, g\bigg\} \bigg/ \bigg\{(n \, x ,u(x))~~\vert~~ x\in X_i\bigg\}.
\end{eqnarray}

\begin{prop}\label{l-adic}
To each biextension $ ({\mathcal{B}}, \Psi_1, \Psi_2,\lambda)$ of $(M_1,M_2)$ by $M_3$ is associated a morphism
\[{\T}_{{\ZZ}/n{\ZZ}}(M_1)\otimes {\T}_{{\ZZ}/n{\ZZ}}(M_2) \longrightarrow   {\T}_{{\ZZ}/n{\ZZ}}(M_3).\]
\end{prop}

\begin{proof} Consider a biextension $ ({\mathcal{B}}, \Psi_1, \Psi_2,\lambda)$ of $(M_1,M_2)$ by $M_3$ and for $i=1,2$ let $m_i$ be an element of ${\T}_{{\ZZ}/n{\ZZ}}(M_i)$ represented by $(x_i,g_i)$ with $u_i(x_i)=n \, g_i.$
 The morphism $\lambda: X_1 \times X_2 \longrightarrow X_3$ gives an element $\lambda(x_1,x_2)$ of $X_3.$ The trivializations $\Psi_1, \Psi_2$ furnish two isomorphisms $a_1$ and $a_2$ from the biextension ${\mathcal{B}}^{\otimes n}$ to the trivial torsor $G_3$:
\begin{eqnarray}
\nonumber  a_1 &:& {\mathcal{B}}^{\otimes n}_{g_1,g_2} \cong {\mathcal{B}} _{ng_1,g_2}={\mathcal{B}}_{u_1(x_1),g_2}
{\buildrel \cong \over \longrightarrow}  G_3 \\
\nonumber  a_2 &:& {\mathcal{B}}^{\otimes n}_{g_1,g_2} \cong {\mathcal{B}} _{g_1,ng_2}={\mathcal{B}}_{g_1,u_2(x_2)}
{\buildrel \cong \over \longrightarrow}  G_3.
\end{eqnarray}
Let
\begin{equation}\label{lh}
a_2 = \phi(m_1,m_2) + a_1 .
\end{equation}
The element $\phi(m_1,m_2)$ of $G_3$ doesn't depend on the choice of the $(x_i,g_i)$ for $i=1,2$.
Because of the compatibility of $u_3 \circ \lambda$ with the trivialization
$(u_1,id_{X_2})^*\Psi_2=(id_{X_1},u_2)^*\Psi_1,$ we
observe that $u_3( \lambda(x_1,x_2))=n \, \phi(m_1,m_2)$. Starting from the biextension $({\mathcal{B}}, \Psi_1, \Psi_2,\lambda)$ we have therefore defined a morphism
\begin{eqnarray}
\nonumber \Phi:{\T}_{{\ZZ}/n{\ZZ}}(M_1)\otimes {\T}_{{\ZZ}/n{\ZZ}}(M_2)  & \longrightarrow &  {\T}_{{\ZZ}/n{\ZZ}}(M_3)\\
\nonumber (m_1,m_2)  & \longmapsto & (\lambda(x_1,x_2), \phi(m_1,m_2)).
\end{eqnarray}
\end{proof}

Recall that the $\ell$-adic realization ${\T}_{\ell}(M_i)$ of the 1-motive $M_i$ is the projective limit of the ${\ZZ}/\ell^n{\ZZ}\,$-modules ${\T}_{{\ZZ}/\ell^n{\ZZ}}(M_i)$ (\cite{D1} (10.1.5)). Using the above proposition, to each biextension of $(M_1,M_2)$ by $M_3$ is associated a morphism $ {\T}_{\ell}(M_1)\otimes {\T}_{\ell}(M_2) \rightarrow {\T}_{\ell}(M_3)$ from the tensor product of the $\ell$-adic realizations of $M_1$ and $M_2$ to the $\ell$-adic realization of $M_3.$

\subsection{Construction of the De Rham realization of biextensions}

Let $S$ be the spectrum of a field $k$ of characteristic 0 embeddable in $\CC$. Let $G_i$ (for $i=1,2$) be a smooth commutative $k$-algebraic group. Let $E$ be an extension of $G_1$ by $G_2$. We can see it as a $G_2$-torsor over $G_1$ endowed with an isomorphism $\nu: pr_1^*E + pr_2^*E \rightarrow \mu^*E$ of $G_2$-torsors over $G_1 \times G_1$, where
$\mu:G_1 \times G_1 \rightarrow G_1$ is the group law on $G_1$
and $pr_i:G_1 \times G_1 \rightarrow G_1$ is the projection (for $i=1,2)$. A \textbf{$\natural$-structure on the extension} $E$ is a connection $\Gamma$ on the $G_2$-torsor $E$ over $G_1$ such that the application $\nu$ is horizontal, i.e. such that $\Gamma$ and $\nu$ are compatible. A \textbf{$\natural$-extension} $(E, \Gamma)$ is an extension endowed with a $\natural$-structure.

Let $G_i$ (for $i=1,2,3$) be a smooth commutative $k$-algebraic group. Let ${\mathcal{P}}$ be a biextension of $(G_1,G_2)$ by $G_3$. We can see it as
 $G_3$-torsor over $G_1 \times G_2$ endowed with an isomorphism $\nu_1: pr_{13}^*{\mathcal{P}} + pr_{23}^*{\mathcal{P}} \rightarrow (\mu_1 \times Id)^*{\mathcal{P}}$ of $G_3$-torsors over $G_1 \times G_1 \times G_2$ and an isomorphism $\nu_2: pr_{12}^*{\mathcal{P}} + pr_{13}^*{\mathcal{P}} \rightarrow (Id \times \mu_2)^*{\mathcal{P}}$ of $G_3$-torsors over $G_1 \times G_2 \times G_2$, which are compatible is the sense of~\cite{SGA7} Expos\'e VII (2.1.1) (here $\mu_i:G_i \times G_i \rightarrow G_i$ is the group law on $G_i$ (for $i=1,2)$, $pr_{i3}:G_1 \times G_1 \times G_2 \rightarrow G_1 \times G_2$ are the projections on the first and second factor for $i=1,2$
and $pr_{1j}:G_1 \times G_2 \times G_2 \rightarrow G_1 \times G_2$
are the projections on the second and third factor for $i=2,3$).
A \textbf{$\natural$-1-structure} (resp. a \textbf{$\natural$-2-structure})
\textbf{on the biextension} ${\mathcal{P}}$ is a connection on the $G_3$-torseur ${\mathcal{P}}$ over $G_1 \times G_2$ relative to $G_1 \times G_2 \rightarrow G_2$ (resp. $G_1 \times G_2 \rightarrow G_1$), such that the applications $\nu_1$ and $\nu_2$ are horizontal.
A \textbf{$\natural$-structure on the biextension ${\mathcal{P}}$} is a
$\natural$-1-structure and a $\natural$-2-structure on ${\mathcal{P}}$, i.e.
a connection $\Gamma$ on the $G_3$-torsor ${\mathcal{P}}$ over $G_1 \times G_2$
such that the applications $\nu_1$ and $\nu_2$ are horizontal.
A \textbf{$\natural$-biextension} $({\mathcal{P}}, \Gamma)$ is an biextension endowed with a $\natural$-structure. The \textbf{curvature $R$ of a $\natural$-biextension $({\mathcal{P}}, \Gamma)$} is the curvature of the underlying connection $\Gamma$: it is a 2-form over $G_1 \times G_2$ invariant by translation  and with values in ${{\li}}(G_3)$, i.e. an alternating form
\[ R:\bigg({\li}(G_1) \times {\li}(G_2)\bigg) \times \bigg({\li}(G_1) \times {\li}(G_2)\bigg) \longrightarrow {\li}(G_3).\]
Since the curvature of the connection underlying a $\natural$-extension
is automatically trivial, the restriction of $R$ to ${\li}(G_1)$ and to ${\li}(G_2)$ is trivial and therefore $R$ defines a pairing (called again \textbf{``the curvature of $({\mathcal{P}},\Gamma)$''})
\begin{equation}\label{dR1}
\Upsilon: {\li}(G_1)\otimes {\li}(G_2) \longrightarrow  {\li}(G_3)
\end{equation}
with $R(g_1+g_2,g'_1+g'_2)=\Upsilon(g_1,g'_2)-\Upsilon(g'_1,g_2).$

Let $K_i=[A_i \rightarrow B_i]$ (for $i=1,2,3$) be a complex of smooth commutative groups. A \textbf{$\natural$-biextension of $(K_1,K_2)$ by $K_3$}
is a biextension of $(K_1,K_2)$ by $K_3$ (see definition \ref{defbiext}) such that the underlying biextension  of $(B_1,B_2)$ by $B_3$ is equipped with a $\natural$-structure and the underlying trivializations are trivializations of $\natural$-biextensions.
 The \textbf{curvature $R$ of a $\natural$-biextension of $(K_1,K_2)$ by $K_3$} is the curvature of the underlying $\natural$-biextension of $(B_1,B_2)$ by $B_3$, or the pairing $\Upsilon: {\li}(B_1)\otimes {\li}(B_2) \rightarrow  {\li}(B_3)$ defined by it.

\begin{lem}\label{dR2}
Let $G_i$ (for $i=1,2$) be an extension of an abelian $k$-variety $A_i$ by a $k$-torus $T_i$. Each extension of $G_1$ by $G_2$ admits a $\natural$-structure.
\end{lem}

\begin{proof}  From the exact sequence
 $0 \rightarrow T_i \rightarrow G_i \rightarrow A_i \rightarrow 0$,
 we have the long exact sequences
\[\begin{array}{c}
 0 \rightarrow {\Hom} (G_1,T_2) \rightarrow {\Hom} (G_1,G_2)
\rightarrow {\Hom} (G_1,A_2) \rightarrow \\
\rightarrow {\Ext}^1 (G_1,T_2)\rightarrow {\Ext}^1 (G_1,G_2)
\rightarrow {\Ext}^1(G_1,A_2) \\
\end{array}\]

\[\begin{array}{c}
 0 \rightarrow {\Hom} (A_1,T_2) \rightarrow {\Hom} (G_1,T_2)
\rightarrow {\Hom} (T_1,T_2) \rightarrow \\
\rightarrow {\Ext}^1 (A_1,T_2)\rightarrow {\Ext}^1 (G_1,T_2)
\rightarrow {\Ext}^1(T_1,T_2) \\
\end{array}\]
According~\cite{S} 7.4 Corollary 1, the group of extensions of $G_1$ by $A_2$  is a torsion group and so modulo torsion, from the first long exact sequence we have the surjection ${\Ext}^1 (G_1,T_2) \rightarrow {\Ext}^1 (G_1,G_2).$ Since ${\Ext}^1(T_1,T_2)=0$, from the second long exact sequence we get a second surjection
${\Ext}^1 (A_1,T_2)\rightarrow {\Ext}^1 (G_1,T_2).$ Therefore after the multiplication by an adequate integer,
each extension of $G_1$ by $G_2$
comes from an extension of the underlying abelian variety $A_1$ by the underlying torus $T_2$. Since the multiplication by an integer for extensions can be viewed as a push-down or a pull-back, and since each extension
of an abelian variety by a torus admits a $\natural$-structure, by pull-back and push-down we get a $\natural$-structure
on each extension of $G_1$ by $G_2$.
\end{proof}

Let $M_i=[X_i {\buildrel u_i \over \longrightarrow} G_i]$ (for $i=1,2,3$)
 be a 1-motive over $k$. The De Rham realization ${\T}_{\rm dR}(M_i)$ of $M_i$ is the Lie algebra of $G^\natural_i$ where $M^\natural_i=[X_i \rightarrow G^\natural_i]$ is the universal vectorial extension of $M_i$ (see~\cite{D1} (10.1.7)). The Hodge filtration on ${\T}_{\rm dR}(M_i)$ is defined by ${\F}^0{\T}_{\rm dR}(M_i)= \ker ( {\li}G_i^\natural \rightarrow {\li}G_i ).$

\begin{prop}\label{naturalstruct}
Each biextension $({\mathcal{B}}, \Psi_1, \Psi_2,\lambda)$ of $(M_1,M_2)$ by $M_3$ defines a
biextension $({\mathcal{B}}^\natural, \Psi_1^\natural, \Psi_2^\natural,\lambda^\natural)$ of $(M_1^\natural,M_2^\natural)$ by $M_3^\natural$ which is endowed with a unique $\natural$-structure.
\end{prop}

\begin{proof} Let $({\mathcal{B}}, \Psi_1, \Psi_2,\lambda)$ be a biextension
of $(M_1,M_2)$ by $M_3$. The proof of this proposition consists of several steps.\\
 \textbf{Step 1:}
Proceeding as in~\cite{D1} (10.2.7.4),
in this step we construct a $\natural$-structure
on the biextension of $(M_1^\natural,M_2^\natural)$ by $M_3$ which is the pull-back of the biextension $({\mathcal{B}}, \Psi_1, \Psi_2,\lambda)$
via the structural projection $M_i^\natural \rightarrow M_i$ for $i=1,2$. For each point $g_1$ of $ G_1$,
${\mathcal{B}}_{g_1}$ is an extension of $G_2$ by $G_3$, which admits a $\natural$-structure according to Lemma \ref{dR2}.
Let $C_{g_1}$ be the set of
$\natural$-structures on ${\mathcal{B}}_{g_1}$. Since two $\natural$-structures differ by an invariant form on $G_2$, $C_{g_1}$ is a torsor under ${\li}(G_2)^*$. The sets $\{C_{g_1}\}_{g_1 \in G_1}$ are the fibers of a ${\li}(G_2)^*$-torsor $C$ on $G_1$. Moreover the Baer'sum of $\natural$-extensions endowed $C$ with a structure of extension of $G_1$ by ${\li}(G_2)^*.$ We lift the morphism $u_1: X_1 \rightarrow G_1$ to $u_1':X_1
\rightarrow C$ in the following way: to each $x_1 \in X_1$ we associate the trivial connection of the trivialized extension ${\mathcal{B}}_{u_1(x_1)}.$
Hence we get a $\natural$-2-structure on the biextension of $([u_i': X_1 \rightarrow C],M_2)$ by $M_3$ pull-back of
$({\mathcal{B}}, \Psi_1, \Psi_2,\lambda)$ via $[u_i': X_1 \rightarrow C] \rightarrow M_1$. By the universal property of $M_1^\natural$ we have a unique commutative diagram
\[\begin{array}{cccccccc}
 0 & \rightarrow & F^0{\T}_{\rm dR}(M_1) & \rightarrow & M_1^\natural & \rightarrow & M_1 & \rightarrow 0   \\
 & & \downarrow  &  & {\scriptstyle v} \downarrow & & = & \\
0 & \rightarrow & {\li}(G_2)^*  & \rightarrow & [u_i': X_1 \longrightarrow C] & \rightarrow & M_1 & \rightarrow 0 . \\
\end{array}\]
If we take the pull-back via $v,$ we have a $\natural$-2-structure on the biextension of $(M_1^\natural,M_2)$ by $M_3.$ Then taking the pull-back via the structural projection $M_2^\natural \rightarrow M_2$ we obtain finally a $\natural$-2-structure on the biextension of $(M_1^\natural,M_2^\natural)$ by $M_3$. Symmetrically
we get a $\natural$-1-structure on this biextension of $(M_1^\natural,M_2^\natural)$ by $M_3$ and hence a $\natural$-structure. \\
\textbf{Step 2:} In this step we show that any biextension
of $(M_1^\natural,M_2^\natural)$ by $M_3$ is canonically the push-down via $M_3^\natural \longrightarrow M_3$ of
a biextension of $(M_1^\natural,M_2^\natural)$ by $M_3^\natural$. In this way  we get a $\natural$-structure on the biextension of $(M_1^\natural,M_2^\natural)$ by $M_3^\natural$ whose push-down is the biextension of $(M_1^\natural,M_2^\natural)$ by $M_3$ of step 1 coming from $({\mathcal{B}}, \Psi_1, \Psi_2,\lambda)$. By definition of the de Rham realization, for $i=1,2,3$ we have the following diagram
\[\begin{array}{cccccccc}
& & 0  &  & 0 & & 0 & \\
& & \uparrow  &  &  \uparrow & &\uparrow & \\
0& \rightarrow & G_i & \rightarrow &  M_i & \rightarrow & X_i &\rightarrow 0 \\
& & \uparrow  &  &  \uparrow & & \cong & \\
0& \rightarrow& G_i' & \rightarrow &  M_i^\natural & \rightarrow & X_i \otimes k& \rightarrow 0 \\
& & \uparrow  &  &  \uparrow & & = & \\
 0 & \rightarrow & F^0 \cap {\W}_{-1} {\T}_{\rm dR}(M_i) & \rightarrow &
F^0 {\T}_{\rm dR}(M_i) & \rightarrow &F^0 \cap {\Gr}_{0} {\T}_{\rm dR}(M_i) & \rightarrow 0   \\
 & & \uparrow  &  &  \uparrow & &\uparrow & \\
& & 0  &  & 0 & & 0 & \\
\end{array}\]
where $G_i'$ is the universal extension of $G_i$. Since
${\Gr}_{0} {\T}_{\rm dR}(M_i)= {\Gr}_{0} {\T}_{\rm dR}(M_i^\natural),$
in order to show that each biextension of $(M_1^\natural,M_2^\natural)$ by $M_3$ lifts to a biextension of $(M_1^\natural,M_2^\natural)$ by $M_3^\natural$ we can restrict to the step ${\W}_{-1},$ i.e. to prove that ${\bBiext} (G'_1,G_2';G_3^\natural)\cong
{\bBiext}(G'_1,G_2';G_3).$
From the short exact sequences
$0 \rightarrow F^0  {\T}_{\rm dR}(M_3) \rightarrow G_3^\natural \rightarrow G_3 \rightarrow 0$, we get the long exact sequences
\[\begin{array}{c}
 0 \rightarrow {\Biext}^0 (G'_1,G_2';F^0  {\T}_{\rm dR}(M_3))\rightarrow {\Biext}^0 (G'_1,G_2';G_3^\natural)\rightarrow {\Biext}^0 (G'_1,G_2';G_3) \rightarrow \\
\rightarrow {\Biext}^1 (G'_1,G_2';F^0 {\T}_{\rm dR}(M_3))\rightarrow {\Biext}^1 (G'_1,G_2';G_3^\natural)\rightarrow
{\Biext}^1(G'_1,G_2';G_3)\rightarrow \\
 \rightarrow {\Ext}^2(G_1' {\buildrel {\scriptscriptstyle \LL}
 \over \otimes}G_2',F^0 {\T}_{\rm dR}(M_3))) \rightarrow ...\\
\end{array}\]
Since for $j=0,1$ and $i=1,2,$ ${\Ext}^j(G_i',{\GG}_a)=0$, we have that
\begin{eqnarray}
\nonumber   {\Biext}^0 (G'_1,G_2';F^0  {\T}_{\rm dR}(M_3)) &\cong&  {\Hom}(G'_1\otimes G_2',F^0{\T}_{\rm dR}(M_3))=0, \\
\nonumber  {\Biext}^1(G'_1,G_2';F^0  {\T}_{\rm dR}(M_3)) &\cong& {\Hom}(G'_1,{\uExt}^1(G'_2,F^0  {\T}_{\rm dR}(M_3)))=0, \\
\nonumber {\Ext}^2(G_1' {\buildrel {\scriptscriptstyle \LL}
 \over \otimes}G_2',F^0 {\T}_{\rm dR}(M_3))) &\cong& {\Ext}^2(G_1', {\RR}{\uHom}( G_2',F^0 {\T}_{\rm dR}(M_3)))=0,
\end{eqnarray}
(the second equivalence is due to~\cite{SGA7} Expos\'e VIII 1.4 and the third one is due to the Cartan isomorphism). Therefore ${\bBiext} (G'_1,G_2';G_3^\natural)\cong
{\bBiext}(G'_1,G_2';G_3).$\\
\textbf{Step 3:} To prove the uniqueness of the $\natural$-structure of the biextension of $(M_1^\natural,M_2^\natural)$ by $M_3^\natural$ coming from $({\mathcal{B}}, \Psi_1, \Psi_2,\lambda)$,
 it is enough to show that any $\natural$-structure on the trivial biextension of $(M_1^\natural,M_2^\natural)$ by $M_3^\natural$ is trivial.
Since the proof is very similar to the one given in~\cite{D1} (10.2.7.4) in the case of biextensions of $(M_1^\natural,M_2^\natural)$ by ${\GG}_m$, we don't give it.
\end{proof}

\begin{cor}\label{dr}
To each biextension $ ({\mathcal{B}}, \Psi_1, \Psi_2,\lambda)$ of $(M_1,M_2)$ by $M_3$ is associated a morphism
\[{\T}_{\dR}(M_1)\otimes {\T}_{\dR}(M_2) \longrightarrow  {\T}_{\dR}(M_3).\]
Explicitly, this morphism is the opposite of the curvature
$\Upsilon: {\li}(G^\natural_1)\otimes {\li}(G^\natural_2) \rightarrow  {\li}(G^\natural_3)$ of the $\natural$-biextension
of $(M_1^\natural,M_2^\natural)$ by $M_3^\natural$ induced by $({\mathcal{B}}, \Psi_1, \Psi_2,\lambda)$.
\end{cor}

\subsection{Comparison isomorphisms}

\begin{prop}\label{compaiso}
(1) Over ${\CC}$, the morphism \ref{l-adic} can be recovered from the morphism \ref{hodge} by reduction modulo $n$.\\
(2) Over ${\CC}$, the morphism \ref{dr} is the complexified of the morphism \ref{hodge}.
\end{prop}

\begin{proof}(1) Recall that by \cite{D1} (10.1.6.2),
 ${\T}_{\ZZ}(M_i) / n {\T}_{\ZZ}(M_i) \cong {\T}_{{\ZZ}/n{\ZZ}}(M_i)$  for $i=1,2,3$. So the assertion follows from the confrontation of (\ref{hl}) and (\ref{lh}).\\
(2) We proceed as in~\cite{D1} (10.2.8). Let $({\mathcal{B}}, \Psi_1, \Psi_2,\lambda)$ be a biextension of $(M_1,M_2)$ by $M_3$. It defines a $\natural$-biextension $({\mathcal{B}}^{\natural}, \Psi_1^{\natural}, \Psi_2^{\natural},\lambda^{\natural})$ of $(M_1^\natural,M_2^\natural)$ by $M_3^\natural$ (see Proposition~\ref{naturalstruct}) and a biextension  $({\mathcal{B}}, \Psi_1, \Psi_2,\lambda)^{\an}$ of $([ {\T}_{\ZZ}(M_1) \oplus {\F}^0{\T}_{\CC}(M_1) \rightarrow {\T}_{\CC}(M_1)],$\\
$ [{\T}_{\ZZ}(M_2) \oplus {\F}^0{\T}_{\CC}(M_2) \rightarrow {\T}_{\CC}(M_2)])$ by
 $[ {\T}_{\ZZ}(M_3) \oplus
 {\F}^0{\T}_{\CC}(M_3) \rightarrow {\T}_{\CC}(M_3)]$ (see Proposition~ \ref{We}). In the analytic category, the $\natural$-biextension $({\mathcal{B}}^{\natural},\Psi_1^{\natural}, \Psi_2^{\natural},\lambda^{\natural})$ defines a $\natural$-structure on the biextension $({\mathcal{B}}, \Psi_1, \Psi_2,\lambda)^{\an}.$ A $\natural$-structure on $({\mathcal{B}}, \Psi_1, \Psi_2,\lambda)$ is in particular a connection on the trivial biextension ${\mathcal{B}}$ of $({\T}_{\CC}(M_1), {\T}_{\CC}(M_2))$ by ${\T}_{\CC}(M_3)$
 i.e. a field of forms $\Gamma_{t_1,t_2}(t'_1+t'_2)$ on the trivial
${\T}_{\CC}(M_3)$-torseur over ${\T}_{\CC}(M_1) \times {\T}_{\CC}(M_2).$
By definition, this connection defines a $\natural$-structure on ${\mathcal{B}}$
if and only if it is compatible with the two group laws underlying ${\mathcal{B}}$, i.e. if and only if
\[ \Gamma_{t_1,t_2}(t'_1+t'_2)= \gamma_1(t_1,t'_2) + \gamma_2(t'_1,t_2) \]
with $\gamma_1$ and $\gamma_2$ bilinear. Moreover in order to have a
$\natural$-structure on $({\mathcal{B}}, \Psi_1, \Psi_2,\lambda)^{\an}$ we have to require that the connection $\Gamma_{t_1,t_2}(t'_1+t'_2)$ and the
trivializations $\Psi_1$ and $\Psi_2$ are compatible, i.e.
$\Psi_1$ and $\Psi_2$ have to be horizontal, and this happens if and only if
$\gamma_i=-\Psi_i$ for $i=1,2$. The curvature of the connection $\Gamma_{t_1,t_2}(t'_1+t'_2)$ is the field of 2-forms:
\begin{eqnarray}
\nonumber  d\Gamma=R_{t_1,t_2}(t_1'+t_2',t_1''+t_2'')  &=&
\Gamma_{t_1',t_2'}(t''_1+t''_2)-\Gamma_{t''_1,t''_2}(t'_1+t'_2)\\
\nonumber    &=&\gamma_1(t'_1,t''_2) + \gamma_2(t''_1,t'_2)-[\gamma_1(t''_1,t'_2) + \gamma_2(t'_1,t''_2)]\\
\nonumber    &=&\gamma_1(t'_1,t''_2)-\gamma_2(t'_1,t''_2)-\big[
\gamma_1(t''_1,t'_2) -\gamma_2(t''_1,t'_2) \big]
\end{eqnarray}
\pn Hence the curvature $\Upsilon:{\T}_{\CC}(M_1) \otimes {\T}_{\CC}(M_2) \longrightarrow  {\T}_{\CC}(M_3)$ of the $\natural$-biextension \\
$({\mathcal{B}}, \Psi_1, \Psi_2,\lambda)^{\an}$ (see (\ref{dR1})) is the form
\[ \Upsilon(t'_1,t''_2)=\gamma_1(t'_1,t''_2) - \gamma_2(t'_1,t''_2) = -\big(\phi_1(t'_1,t''_2)-\phi_2(t'_1,t''_2)\big)=-\Phi(t'_1,t''_2).\]
\end{proof}

\subsection{Compatibility with the category of mixed realizations}

Let $S$ be the spectrum of a field $k$
of characteristic 0 embeddable in $\CC.$ Fix an algebraic closure $\ok$ of $k$. Let ${\mathcal{MR}}_{\ZZ}(k)$ be the integral version of the neutral Tannakian category over
$\QQ$ of mixed realizations
(for absolute Hodge cycles) over $k$ defined by Jannsen in~\cite{J} I 2.1.
The objects of ${\mathcal{MR}}_{\ZZ}(k)$ are families
\[N=((N_\sigma, \mathcal{L}_\sigma),N_{\rm dR},N_{\ell}, I_{\sigma, {\rm dR}}, I_{{\overline \sigma}, {\ell}} )_{\ell,\sigma,{\overline \sigma}}\]
where
\begin{itemize}
  \item $N_\sigma$ is a mixed Hodge structure for any embedding
 $\sigma:k \rightarrow \CC$ of $k$ in $\CC$;
  \item $N_{\rm dR}$ is a finite dimensional $k$-vector space with an increasing filtration ${\W}_*$ (the Weight filtration) and a decreasing filtration ${\F}^*$ (the Hodge filtration);
  \item $N_{\ell}$ is a finite-dimensional $\QQ_\ell$-vector space with a continuous $\gal$-action and an increasing filtration ${\W}_*$ (the Weight filtration), which is $\gal$-equivariant, for any prime number $\ell$;
  \item $I_{\sigma, {\rm dR}}:N_\sigma \otimes_{\QQ} {\CC} \rightarrow
N_{\rm dR}\otimes_k {\CC} $
and $I_{{\overline \sigma}, {\ell}}:N_\sigma \otimes_{\QQ}
{\QQ}_\ell \rightarrow N_{\ell}$ are comparison isomorphisms for any $\ell$, any $\sigma$ and any
${\overline \sigma}$ extension of $\sigma$ to the algebraic closure of $k$;
  \item $\mathcal{L}_\sigma$ is a lattice in $ N_\sigma$ such that,
for any prime number $\ell$, the image $\mathcal{L}_\sigma \otimes {\ZZ}_\ell$
of this lattice through the comparison isomorphism $I_{{\overline \sigma}, {\ell}}$ is a $\gal$-invariant subgroup of $N_{\ell}$ ($\mathcal{L}_\sigma$ is the integral structure of the object $N$ of ${\mathcal{MR}}_{\ZZ}(k)$).
\end{itemize}
Before to define the morphisms of the category ${\mathcal{MR}}_{\ZZ}(k)$ we have to introduce the notion of Hodge cycles and of absolute Hodge cycles. Let $N=((N_\sigma, \mathcal{L}_\sigma),N_{\rm dR},N_{\ell}, $\\
$I_{\sigma, {\rm dR}}, I_{{\overline \sigma}, {\ell}} )_{\ell,\sigma,{\overline \sigma}}$ be an object of the Tannakian category
 ${\mathcal{MR}}_{\ZZ}(k)$.
A Hodge cycle of $N$ relative to an embedding $\sigma : k \rightarrow {\CC}$  is an element $(x_\sigma,x_{\rm dR},x_{\ell})_{\ell}$ of $N_\sigma \times N_{\rm dR} \times \prod_{\ell} N_{\ell}$ such that $I_{\sigma, {\rm dR}}(x_\sigma) =x_{\rm dR}, I_{{\overline \sigma}, {\ell}}(x_\sigma) =x_{\ell}$ for any prime number $\ell$ and $x_{\rm dR} \in {\F}^0 N_{\rm dR} \bigcap {\W}_0 N_{\rm dR}. $
An absolute Hodge cycle is a Hodge cycle relative to every embedding $\sigma : k \rightarrow {\CC}$. By definition, the morphisms of the Tannakian category ${\mathcal{MR}}_{\ZZ}(k)$ are the absolute Hodge cycles: more precisely, if $N$ and $N'$ are two objects of ${\mathcal{MR}}_{\ZZ}(k)$, the morphisms ${\Hom}_{{\mathcal{MR}}_{\ZZ}(k)}(H,H')$ are the absolute Hodge cycles of the object ${\uHom}(H,H')$ (see~\cite{J} I Definition 2.1 and (2.11)).

Since 1-motives are endowed with an integral structure, according to~\cite{D1} (10.1.3) we have the fully faithful functor
\begin{eqnarray}
\nonumber   \{\mathrm{1-motives~}/~k\} & \longrightarrow & {\mathcal{MR}}_{\ZZ}(k) \\
\nonumber  M & \longmapsto &
{\T}(M)=(({\T}_\sigma(M),\mathcal{L}_\sigma),{\T}_{\rm dR}(M),{\T}_{\ell}(M),
I_{\sigma, {\rm dR}}, I_{{\overline \sigma}, {\ell}} )_{\ell, \sigma,
  {\overline \sigma}}
\end{eqnarray}
\pn which attaches to each 1-motive $M$ of ${\mathcal{M}}(k)$
its Hodge realization $({\T}_\sigma(M),\mathcal{L}_\sigma)$ with integral structure for any embedding
 $\sigma:k \rightarrow \CC$ of $k$ in $\CC$, its de Rham realization ${\rm T}_{\rm dR}(M)$, its $\mathbf{\ell}$-adic realization ${\rm T}_{\ell}(M)$ for any prime number $\ell$, and its
comparison isomorphisms.

\begin{thm}\label{realthm}
Let $M_i$ ($i=1,2,3$) be a 1-motive over $k$. We have that
\[{\Hom}(M_1, M_2;M_3)\otimes {\QQ} \cong {\Hom}_{{\mathcal{MR}}_{\ZZ}(k)}\big({\T}(M_1)\otimes {\T}(M_2), {\T}(M_3)\big)\]
\end{thm}

\begin{proof} Let $(M_i)_\sigma = M_i \otimes_{\sigma} {\CC}$ (for $i=1,2,3$). According to Corollary~\ref{hodge}, the biextensions of $((M_1)_{\sigma},(M_2)_{\sigma})$ by $(M_3)_{\sigma}$ are bilinear morphisms ${\T}_{\sigma}(M_1)\otimes {\T}_{\sigma}(M_2) \rightarrow  {\T}_{\sigma}(M_3)$ in the category $\mathcal{MHS}$ of mixed Hodge structures, i.e. they are rational tensors living in
$$ {\F}^0 \bigcap {\W}_0 \bigg({\T}_{\sigma}(M_1)\otimes {\T}_{\sigma}(M_2) \otimes {\uHom}_{\mathcal{MHS}}({\T}_{\sigma}(M_3), {\T}_{\sigma}({\ZZ}))\bigg).$$
 By Proposition~\ref{l-adic} and Corollary~\ref{dr}, biextensions define bilinear morphisms also in the $\ell$-adic and De Rham realizations and all these bilinear morphisms are compatible through the comparison isomorphisms (see Proposition~\ref{compaiso}). Therefore biextensions of 1-motives define Hodge cycles. In~\cite{Br} Theorem (2.2.5) Brylinski proves that Hodge cycles over a 1-motive defined over $\ok$ are absolute Hodge cycles and so biextensions are Hodge cycles relative to every embedding ${\overline \sigma}: \ok \rightarrow {\CC}$. Since biextensions of $(M_1,M_2)$ by $M_3$ are defined over $k$, the bilinear morphisms they define are invariant under the action of $\gal$. This implies that biextensions $(M_1,M_2)$ by $M_3$ are Hodge cycles relative to every embedding $ \sigma: k \rightarrow {\CC},$ i.e. they are
morphisms in the category ${\mathcal{MR}}_{\ZZ}(k)$.
\end{proof}



\end{document}